\documentclass[numbers,preprint,10pt]{elsarticle}
\usepackage{amsmath,amsfonts}
\DeclareMathOperator*{\argmax}{\arg\!\max}
\DeclareMathOperator*{\argmin}{\arg\,\min}
\usepackage{xcolor}
\usepackage{pgfgantt}
\usepackage[font=small]{subcaption}
\usepackage{hyperref}
\usepackage{lscape}
\usepackage{float}
\usepackage{graphicx}

\biboptions{numbers}

\hypersetup{colorlinks,
citecolor=black
}

\usepackage{tikz}
\usetikzlibrary{positioning}
\usetikzlibrary{external}
\usetikzlibrary{patterns}
\definecolor{darkgreen}{RGB}{0, 204, 0}
\definecolor{darkred}{RGB}{214, 41, 41}
\definecolor{darkblue}{RGB}{34, 81, 221}

\usepackage{bm}

\usetikzlibrary{positioning}

\definecolor{cardiffRed}{cmyk}{0.0,0.91,0.65,0.11}
\definecolor{AMCG_blue}{RGB}{19,12,168}
\definecolor{MUFFINS_blue}{RGB}{0,114,182} % light blue
\definecolor{darkgreen}{RGB}{0, 204, 0}
\definecolor{darkred}{RGB}{214, 41, 41}
\definecolor{darkblue}{RGB}{0,114,182}% 
\definecolor{MUFFINS_yellow}{RGB}{233,143,0}

\usepackage{amsmath,amssymb,amsfonts}
\usepackage{caption}
\usepackage{algorithm}
\usepackage{algorithmicx}
\usepackage{algpseudocode}

\usepackage{listings}%
\usepackage{alltt}%
\usepackage{eqparbox}

\usepackage{setspace}

\usepackage{relsize} %make operators larger 

\usepackage{mathtools}

\usepackage{array}
\newcolumntype{L}[1]{>{\raggedright\arraybackslash}m{#1}}
\newcolumntype{R}[1]{>{\raggedleft\arraybackslash}m{#1}}
\newcolumntype{C}[1]{>{\centering\arraybackslash}m{#1}}

\usepackage{siunitx}

\usepackage{easylist}

\usepackage[margin=2cm]{geometry}

\begin{document}

\begin{frontmatter}

\title{An autoencoder-based reduced-order model for eigenvalue  problems with application to neutron diffusion}

\author[IC]{Toby Phillips}
\author[IC]{Claire E. Heaney\corref{cor2}}
%\corref{cor2}\ead{c.heaney@imperial.ac.uk}
\ead{c.heaney@imperial.ac.uk}
\author[JACOBS]{Paul N. Smith} 
\author[IC]{Christopher C. Pain} 

\cortext[cor2]{Corresponding author.}
\address[IC]{Applied Modelling and Computation Group, Department of Earth Science and \\Engineering, Imperial College 
London, UK}
\address[JACOBS]{Jacobs, Poundbury, Dorset, UK}

\begin{abstract}
Using an autoencoder for dimensionality reduction, this paper presents a novel projection-based reduced-order model for eigenvalue problems.  Reduced-order modelling relies on finding suitable basis functions which define a low-dimensional space in which a high-dimensional system is approximated. Proper orthogonal decomposition (POD) and singular value decomposition (SVD) are often used for this purpose and yield an optimal linear subspace. Autoencoders provide a nonlinear alternative to POD/SVD, that may capture, more efficiently, features or patterns in the high-fidelity model results. 

Reduced-order models based on an autoencoder and a novel hybrid SVD-autoencoder are developed. These methods are compared with the standard POD-Galerkin approach and are applied to two test cases taken from the field of nuclear reactor physics.
\end{abstract}


\end{frontmatter}

\section{Introduction}\label{introduction} 
Reduced-order modelling or model reduction~\cite{schilders2008model} is a numerical technique that has made a significant impact across a broad range of fields, including aerodynamics~\cite{Zahr2015}, haemodynamics~\cite{Ballarin2016}, fracture~\cite{Kerfriden2012,Agathos2020}, porous media~\cite{Jansen2017,Xiao2018} and molecular dynamics~\cite{Hoang2016}. Its aim is to produce a low-dimensional model which is a good approximation of a high-dimensional model (high-fidelity model), but which can be solved in a fraction of the computational time required by the high-fidelity model (HFM). Main applications for reduced-order modelling are multi-query problems, such as optimisation and control~\cite{Chinesta2013,Zahr2015}, and applications which require real-time solutions~\cite{Niroomandi2011,Aguado2015}. The computational efficiency of reduced-order models is achieved by an offline/online split of the computational work. There are typically three parts to the offline stage. First, the high-fidelity model is solved for different parameter values. Second, dimensionality reduction techniques are applied to the HFM solutions (referred to as snapshots) producing basis functions that span the low-dimensional space or reduced space. Finally, the discretised governing equations are approximated in the reduced space. In the online stage, the reduced-order model is solved for an unseen parameter value in a relatively fast time. 

% balanced truncation, DMD, Krylov subspace, reduced basis
A number of dimensionality reduction methods have been proposed for finding the basis of the reduced space, however Proper Orthogonal Decomposition~\cite{Taira2017,brunton2020machine} (POD) is by far the most common. Also known as Principal Component Analysis (PCA), the Karhunen-Lo{\`e}ve transform, Empirical Orthogonal Functions and the Hotelling transform, POD was introduced by Lumley~\cite{Lumley1967} as a method to extract coherent structures from turbulent flow data, and was made tractable with the method of snapshots by Sirovich~\cite{sirovich1987turbulence}. The POD basis functions are found by taking a Singular Value Decomposition (SVD) of the snapshots, or by solving a related eigenvalue problem. Using the basis functions, a Galerkin projection can be applied to the discretised governing equations (i.e.~the high-fidelity model) to produce the reduced system of equations. For HFMs whose operators can be affinely decomposed, the Galerkin projection will be applied to all the matrices in the decomposition. Finding the reduced system for an unseen parameter can be done by interpolating the resulting reduced system matrices. Reduced-order models produced by this method are referred to as projection-based reduced-order models~\cite{benner2015survey}. For non-affine or non-linear systems, this process is more complex, and generally more sampling is required to represent the non-affine nature or non-linearities of the system. Hyper-reduction~\cite{Ryckelynck2005,Fritzen2018} and Discrete Empirical Interpolation~\cite{Chaturantabut2010} are two methods that can be used for this. %

In some circumstances, the projection step is challenging, as the matrices representing the high-fidelity model matrices are not easily available, for example when using licensed software or for complex multi-physics problems. To alleviate this difficulty, an alternative method has been developed and is referred to variously as `non-intrusive reduced-order modelling' (NIROM), `POD with interpolation', `Galerkin-free', `surrogate POD' (see~\cite{Benamara2017} and references within), `data-driven reduced-order modelling' \cite{Kaiser2014,Guo2019,Swischuk2019} and `system/model identification' \cite{Polifke2014,Wang2018}. As in projection-based methods, these non-intrusive methods take a set of high-fidelity model solutions (snapshots) and find a reduced basis, often by using POD. However, to find solutions for unseen parameter values, projection of the high-fidelity system onto the basis functions is replaced by interpolating between snapshots (or snapshots that have themselves been projected onto the basis functions). Although cubic interpolation~\cite{Bui-Thanh2003} or radial basis functions~\cite{Xiao2018} were originally used for the interpolation, NIROM methods have since embraced machine learning, using a variety of neural networks for this purpose, including autoencoders in combination with Long Short-Term Memory networks~\cite{Gonzalez2018,Wiewel2019}, multi-layer perceptrons~\cite{Hesthaven2018} and Gaussian Process Regression~\cite{Guo2019,Xiao2019}. 
%, although initially methods used were linear interpolation~\cite{Shinde2016} or cubic spline interpolation~\cite{Bui-Thanh2003,Ly2001}. Radial basis functions have also been used as interpolants~\cite{Xiao2018, Audouze2013, Chen2018}. 
When using neural networks, the offline stage involves training a neural network to learn how the solutions depend on various model parameters, and the online stage involves evaluating the neural network for an unseen parameter value or values.

Although machine learning techniques are becoming more prevalent for predicting the behaviour of the HFM (third part of offline stage and online stage in NIROM-type methods), it is much less common to find machine learning algorithms being used in the dimensionality reduction stage (second part of offline stage in both projection-based or NIROM methods). One type of neural network that is ideal for dimensionality reduction is the autoencoder. A type of unsupervised (or self-supervised) feed-forward neural network, the autoencoder learns the identity map with a network architecture that includes a bottleneck. The central bottleneck layer forces the autoencoder to learn a reduced or low-dimensional representation of the data. The neurons in this layer determine the so-called latent space and are referred to as latent variables. These types of networks are common in image classification and identification~\cite{hinton2006reducing,wang2012folded}, and have been used with great success to fill in gaps in images~\cite{Mao2016} and to remove noise from data~\cite{Pascal2010}. Several authors~\cite{Oja1982,Bourlard1988, Baldi1989} have highlighted the connection between autoencoders and SVD (or PCA), namely that an autoencoder with one hidden layer and \emph{linear} activation functions can produce basis functions than span the same space as basis functions derived from an SVD. 
%\textcolor{red}{Wang 2016 et al}~\cite{Wang2016auto} have also compared the autoencoder with four other dimensionality reduction methods: PCA and linear discriminant analysis (linear methods), and locally linear embedding and the Isomap method (nonlinear methods). \textcolor{red}{conclusion? AE can detect repetitive patterns as well as reduce the dimension} ? \textcolor{red}{Tyler Manning-Dahan}
%The weights associated with the decoder must be constrained to be the transpose of the encoder, and the basis functions will neither be ordered nor orthonormalised.  
This connection between the autoencoder and the SVD paves the way for an autoencoder to be seen as a tool to perform \emph{nonlinear} dimensionality reduction (by choosing nonlinear activation functions) and therefore to be used as an alternative to the SVD/PCA~\cite{Takane1998,hinton2006reducing}. Amongst those early to recognise the potential of the autoencoder for dimensionality reduction within a model reduction framework were Milano~et~al.~\cite{Milano2002}, who reconstructed the near wall field from computational solutions for pressure and shear stress at the wall. 
The suitability of autoencoders for dimensionality reduction having been established, a number of authors have since explored their use in reduced-order models, mostly non-intrusive reduced-order models. This is because the strength of the autoencoder is that it is a nonlinear embedding, but this also means that the reduced-order system will be nonlinear and is therefore more complex to solve than for a POD-based ROM. For non-intrusive models, this nonlinearity does not present any additional challenge and no change in the algorithm is required, hence the greater uptake of autoencoders in NIROM-type methods. For example,  Gonzalez~et~al.~\cite{Gonzalez2018} and Wiewel~et~al.~\cite{Wiewel2019} both use a convolutional autoencoder to reduce the dimension of their problems and Long Short-Term Memory networks to learn the dynamics. The former demonstrate their method on interacting vortices in a 2D box and 2D lid-driven cavity flow. The interacting-vortices problem aims to demonstrate the location invariance properties of convolutional networks. %
The long-term statistics of the flow, as embodied by the turbulent kinetic energy, are captured much better by the autoencoder-based method than by the projection-based POD-Galerkin model. Coming from the field of computer graphics, Wiewel~et~al.~\cite{Wiewel2019} solve 3D problems with millions of spatial degrees of freedom. With their reduced-order model, they are able to produce realistic-looking simulations of complex dynamical systems such as sloshing waves, colliding bodies of fluid and smoke convection. Their reduced model runs approximately two orders of magnitude faster than their high fidelity solver.

As previously mentioned, the substitution of an autoencoder for an SVD is less straightforward for projection-based reduced-order models than it is for non-intrusive reduced-order models. There are a couple of simple examples of the incorporation of autoencoders in a projection-based reduced-order model~\cite{hartman2017deep,kashima2016nonlinear}, however, the most comprehensive to date is a paper by Lee~et~al.~\cite{lee2020model} They propose a framework for projecting dynamical systems onto nonlinear manifolds using deep convolutional autoencoders and derive \emph{a posteriori} error bounds. Focusing on advection-dominated benchmark problems known to be challenging for POD-based ROMs due to a slow decay of singular values or Kolmogorov $n$-width~\cite{Ahmed2020}, they model the convection, diffusion and chemical reactions of a flame as it evolves in time over a range of values for two parameters which occur in the nonlinear reaction source term. The error in the autoencoder-based method is over 40 times less than that of the standard POD-based projection method, showing that, although their method requires a higher offline computational time (for the training of the autoencoder), they can produce much more accurate results than projection-based methods for some nonlinear problems. %These systems are non-linear and highlight the benefit of using autoencoders for this compression.
In this paper, we develop autoencoder-based ROMs for eigenvalue problems, using a standard autoencoder and a hybrid SVD-autoencoder.

%In this paper, we develop autoencoder-based ROMs for eigenvalue problems, using a standard autoencoder, a hybrid SVD-autoencoder and a variational autoencoder (VAE). The VAE is a type of generative model that combines deep learning with statistics. Introduced by Kingma~et~al.~\cite{Kingma2014,Kingma2019}, these networks encode the inputs as a distribution over latent space. It is possible to take a sample from this distribution and to decode the result in order to generate examples that have not been seen in the training process. By contrast with standard autoencoders, the latent space of a VAE is highly structured and smooth, and the latent variables correspond to a meaningful variation in the data~\cite{Chollet}. For example, White~\cite{White2017} trained a VAE with images of celebrities' faces and was able to identify a `smile' vector. As this property was captured by one variable in latent space, he could add and take-away smiles from images. In one of the test cases in this paper, we see whether the VAE is better able to generalise than the standard AE. 
%http://www.dpkingma.com/sgvb_mnist_demo/demo.html

The field of nuclear engineering has some extremely demanding applications such as neutron transport within a reactor core, the geometry of which presents extremely challenging aspect ratios. A recent discretisation of a Westinghouse PWR-900 core with 44 energy groups used over 1.7 trillion degrees of freedom~\cite{palmtag2014coupled,slaybaugh2018eigenvalue} precluding the use of optimisation, safety-analysis or real-time analysis. It comes as no surprise, therefore, that reduced-order modelling has been applied to reactor physics problems. Based on the multigroup neutron diffusion equation for criticality, Chunyu~et~al.~\cite{chunyu2018fast} applied a reduced basis method to 2D and 3D Pressurized Water Reactor~(PWR) benchmark problems. Buchan~et~al.~\cite{buchan2013pod} also tackled criticality by transforming the eigenvalue problem into a time dependent problem and applied this to simple 1D and 2D models of nuclear reactors (test cases based on these are used in this paper). These papers~\cite{chunyu2018fast,buchan2013pod} used a single POD basis for all the energy groups, the so-called monolithic approach. Heaney~et~al.~\cite{Heaney2018} developed a method to model control-rod movement with a POD-based ROM for a PWR fuel assembly, in which POD basis functions were calculated independently for each energy group. German~et~al.~\cite{German2019} demonstrated that such an approach produced more accurate results than the monolithic approach. All the methods discussed in this paragraph show computational gains of at least three orders of magnitude over the respective HFMs. Outside of criticality problems, ROM has also been applied to time-dependent movement of control rods~\cite{sartori2015reduced}, and coupled neutronics and heat transfer~\cite{sartori2016multi}, both using a reduced basis method based on POD.
%, and using POD, Reed~et~al.\cite{reed2019effectiveness} were able to reduce substantially the number of coefficients used to represent the dependence of the neutron flux on energy.

%The reduced basis finite element method is used in~\cite{chunyu2018fast} to solve the multi-group neutron diffusion equation. Sartori \textit{et al.} simulated control rod movement by a ROM based on the reduced basis method~\cite{sartori2016reduced} on a 2D model and advanced this further where a time dependent 3D model of the same reactor is modelled~\cite{sartori2015reduced}. Another paper~\cite{reed2019effectiveness} uses the discrete generalised multi-group method to solve the neutron transport equation. Proper generalized decomposition is used on the multigroup neutron diffusion eigenvalue problem~\cite{prince2020application}. All of these papers reduce the computational cost of their models. In~\cite{lorenzi2018adjoint} POD is seen to perform better than the modal method for reducing the neutron diffusion equation. A ROM based on POD and the method of snapshots was applied to a PWR assembly by Heaney \textit{et al.}~\cite{Heaney2018}. POD with snapshots has been successfully applied to resolving the direction of travel of neutrons and photon problems~\cite{buchan2015pod}. All of these show POD performing exceptionally well on neutron transport problems.

Although, in many areas, adding diffusion or solving diffusion problems often leads to fields that are smooth and well able to be represented accurately with POD basis functions, in other areas, such as reactor physics, one is often confronted by problems with abruptly changing fields similar to advection problems. A good example is a control rod that is partially inserted into a reactor, where one would see a near zero flux within the control rod and a much higher flux a small distance away from the control rod; see the scalar flux solutions shown in Buchan~et~al.~\cite{Buchan2019} for instance. This similarity between advection and diffusion problems is also seen in discretisation methods such as the Self Adjoint Angular Flux method~\cite{Morel1999}, in which the first order transport equations are transformed into a series of diffusion-like equations with the application of a least squares principle and a particular weighting. Although two relatively simple test cases are presented in this paper, our future work will involve increasingly complicated/more realistic assembly or reactor configurations which will benefit from the improved modelling capabilities that autoencoders provide.

This paper reports, what the authors believe to be, the first application of an autoencoder-based reduced-order model to an eigenvalue problem. Also presented is a novel hybrid SVD-autoencoder method that combines an SVD with an autoencoder. This is useful as training a fully-connected autoencoder with a long input vector is computationally expensive. Without the SVD, the length of the input vector would be equal to the number of degrees of freedom of the problem. Applying the SVD reduces the length of the input vector to, at most, the number of snapshots.

The remainder of the paper is organised as follows. In Section~2, the governing equations and their control-volume discretisation are presented. Used to find the dominant eigenvalues and eigenvectors of eigenvalue problems, the power method is also described. Section~3 contains background information relevant to reduced-order modelling techniques, including descriptions of POD, SVD, selecting basis functions and constructing projection-based reduced-order models. A description of autoencoders is included in Section~4 along with theory on how to include these networks in a projection-based ROM. In Section~5 results of the 1D and 2D test cases are given for several different autoencoder-based reduced-order models and they are compared with a standard POD-based reduced-order model. Concluding remarks are made in the final section.

%-----------------------------------------------------------------
%-----------------------------------------------------------------
\section{Governing equations and their discretisation}\label{governing} 
This section introduces the equations that govern criticality and their control-volume discretisation, referred to as the high-fidelity model. A description is then given of the power method, which can be used to solve generalised eigenvalue problems such as that arising in criticality.
\subsection{Diffusion Equation}\label{diffusion_equation} 
The one-group steady-state diffusion equation for criticality takes the form of a generalised eigenvalue problem and can be written as
\begin{equation}\label{eq:diff-eig}
- \nabla\cdot (D \nabla \phi) + \Sigma^a \phi = \lambda \nu \Sigma^f \phi,
\end{equation}
where $\phi$ is the scalar flux of the neutron population, $\Sigma^a$ represents the absorption cross-section, $\Sigma^f$ represents the fission cross-section and $\nu$ is the average number of neutrons produced per fission event. The diffusion coefficient, $D$, is given by:
\begin{equation}\label{D_equation}
    D = \frac{1}{3(\Sigma^a + \Sigma^s)}\,,
\end{equation}
where $\Sigma^s$ is the scattering coefficient. The eigenvalue, $\lambda$, is defined as the reciprocal of $k_{\text{eff}}$, that is $\lambda = \frac{1}{k_{\text{eff}}}$, where:
\begin{equation}
    k_{\text{eff}}=\frac{\text{number of neutrons in one generation}}{\text{number of neutrons in the preceding generation}}.
\end{equation}
The boundary conditions are for reflection: 
\begin{equation}\label{eq:diff-reflect}
D\frac{\partial \phi}{\partial n} = 0\,,
\end{equation}
and for a vacuum or bare-surface:
\begin{equation}\label{eq:diff-bare}
-D\frac{\partial \phi}{\partial n} = \frac{1}{2} \phi\, ,
\end{equation}
in which $n$ is the outward-pointing normal to the boundary. In this paper, autoencoder-based reduced-order models are developed for the one-group equations, however the methods described here could easily be applied to the multi-group equations.

\subsection{Discretisation}
A control-volume discretisation of the diffusion equation in 2D with 
a regular $N_x \times N_y$ mesh can be written as:
% cheating a bit with eqnarray
\begin{eqnarray}%{eq:discretised_hfm}
& &\frac{1}{2} \max\left\{D_{i-1,j} + D_{i,j},\; 0\right\}  \frac{ \phi_{i,j}-\phi_{i-1,j}}{\Delta x^2} \ + \ \frac{1}{2} \max\left\{D_{i,j} + D_{i+1,j},\; 0\right\}  \frac{ \phi_{i+1,j}-\phi_{i,j}}{\Delta x^2} \nonumber\\[2mm] 
 & & + \,\frac{1}{2} \max\left\{D_{i,j-1} + D_{i,j},\; 0\right\}  \frac{ \phi_{i,j}-\phi_{i,j-1}}{\Delta y^2}   
  + \frac{1}{2} \max\left\{D_{i,j} + D_{i,j+1},\; 0\right\}  \frac{ \phi_{i,j+1}-\phi_{i,j}}{\Delta y^2} + 
  {\Sigma^a}_{i,j} \phi_{i,j} \label{eq:discretised_hfm}\\[3mm]
& & = \lambda\nu\Sigma^f_{i,j}\phi_{i,j}, \qquad  \forall i\in\{2,3,..,N_x-1\}, \  \forall j\in\{2,3,..,N_y-1\},\nonumber
\end{eqnarray}
in which $\Delta x$ and $\Delta y$ are the uniform cell widths in the $x$- and $y$-directions respectively, $N_x$ and $N_y$ are the numbers of cells in the $x$- and $y$-directions respectively, the subscripts $i,j$ refer to the cells in the $x$- and $y$-directions respectively, and $\phi_{i,j}$ represents the scalar flux in cell~$i,j$. In this expression, the first and last cells are omitted in order to apply the boundary conditions efficiently. To apply a reflective boundary condition, see Equation~\eqref{eq:diff-reflect}, the diffusion coefficients in the outermost cells are set to a large negative number. To set a bare surface boundary condition, see Equation~\eqref{eq:diff-bare}, the diffusion coefficients in the outermost cells are again set to a large negative number and the absorption term is modified as now described. To apply such a boundary condition to the left or right edges (where the normal to the boundary is aligned with the $x$-direction): 
\begin{equation}\label{eq:abs-bc-x}
\Sigma^a_{i,j}\leftarrow \Sigma^a_{i,j} + \frac{1}{2\Delta x}\,,
\end{equation}
and to apply this to the top or bottom edges (where the normal to the boundary is aligned with the $y$-direction): 
\begin{equation}\label{eq:abs-bc-y}
\Sigma^a_{i,j}\leftarrow \Sigma^a_{i,j} + \frac{1}{2\Delta y}\,.
\end{equation}
For cells that have both boundary conditions the following modification is made: 
\begin{equation}\label{eq:abs-bc-xy}
\Sigma^a_{i,j}\leftarrow \Sigma^a_{i,j} + \frac{1}{2\Delta x} + \frac{1}{2\Delta y}\,. 
\end{equation}

The discretised form of Equation~\eqref{eq:diff-eig} can therefore be written as:
\begin{equation}\label{disc_fiss}
    \bm{A}\bm{\phi}=\lambda \bm{B}\bm{\phi}.
\end{equation}
where matrix~$\bm{A}$ contains the transport terms (scattering, diffusion) from the left-hand side of Equation~\eqref{eq:discretised_hfm},  matrix~$\bm{B}$ represents the fission terms given in the right-hand side of Equation~\eqref{eq:discretised_hfm} and the vector~$\bm{\phi}$ contains the values of the scalar flux for each cell. The matrices are of size $(N_x-2)(N_y-2)$ by $(N_x-2)(N_y-2)$. Although a 2D discretisation is given here, the methods can also be applied in 1D or 3D.

\subsection{The Power Method}\label{sec:power_it}
The power method~\cite{Golub1996} is an algorithm which finds the dominant eigenvalue and eigenvector of an eigenvalue problem. In certain circumstances the method can be slow to converge, however, it is possible to accelerate convergence~\cite{prince2020application}. The power method has a broad range of application and has been used to rank scientific journals~\cite{Pinski1976} and to make suggestions to Twitter users of whom to follow~\cite{Gupta}. In order to apply the power method to the criticality eigenvalue problem~\cite{Scheben}, Equation~\eqref{disc_fiss} is rearranged to be in the form of a standard eigenvalue problem as follows
\begin{equation}\label{eq:rearranged_gevp}
\bm{A}^{-1}\bm{B} \bm{\phi} = \frac{1}{\lambda}\bm{\phi}\,.
\end{equation}
First the algorithm iterates to find a flux solution by repeatedly applying $\bm{A}^{-1}\bm{B}$ to the current scalar flux. Once these inner iterations have converged, the algorithm normalises the scalar flux and calculates the associated eigenvalue. The outer iterations continue until convergence is reached, which, in this paper, is either when the difference in the eigenvalue in successive iterations is less than \num{e-8} or when 1000 iterations have been completed. This procedure is detailed in Algorithm~\ref{alg:PM_outer}, in which $\bm{b}$ is a vector whose entries all equal one. The normalisation applied here, seen on line~\ref{alg:normalisation}, differs from the normalisation usually applied in the power method, however, results in a unit number of fissions, which is often preferred in nuclear applications. The parenthesised superscript on the flux solutions and the eigenvalues indicates the outer iteration index. The inner iteration loop is not described here. Although not described here, the inner iterations use a forward backward Gauss Seidel method to solve Equation~\eqref{eq:rearranged_gevp} when generating the high-fidelity model solutions.

\begin{algorithm}
\caption{Power Method: Outer Iterations}\label{alg:PM_outer}
\begin{algorithmic}[1]
\Function{Outer\textunderscore{}Iterations}{$\bm{A}$, $\bm{B}$, $\bm{\phi}^{\text{guess}}$, $\lambda^{\text{guess}}$} 
\State $\bm{\phi}^{(0)} = \bm{\phi}^{\text{guess}}$, $\lambda^{(0)} = \lambda^{\text{guess}}$
\State $i$\textunderscore{}$max = 1000$
\State $k$\textunderscore{}$tol = 10^{-8}$
\State $i = 0$
\State \texttt{not\textunderscore{}converged = True}
\While{\texttt{not\textunderscore{}converged}}
\setstretch{1.3}
   \State $\bm{\phi}^{(i+1)} = \Call{Inner\textunderscore{}Iterations}{\bm{A}, \bm{B}, \bm{\phi}^{(i)}, \lambda^{(i)}} $ \Comment{to solve Equation~\eqref{eq:rearranged_gevp}}
   \State $\bm{\phi}^{(i+1)} \gets \dfrac{\bm{\phi}^{(i+1)}}{\bm{b}^T \bm{B} \bm{\phi}^{(i+1)}}$ \Comment{normalising the flux$\qquad\qquad\qquad\qquad\qquad\qquad\qquad$} \label{alg:normalisation}
\setstretch{1.6}
   \State $\lambda^{(i+1)} = \dfrac{\bm{b}^T \bm{A} \bm{\phi}^{(i+1)}}{\bm{b}^T \bm{B} \bm{\phi}^{(i+1)}}$ 
\setstretch{2.2} 
   \If{$(i=i$\textunderscore{}$max$ \textbf{or} $|k_{\text{eff}}^{(i+1)} - k_{\text{eff}}^{(i)}| < k$\textunderscore{}$tol)$}
\setstretch{1.6}
     \State \texttt{not\textunderscore{}converged = False}
\setstretch{1.0}
    \State $\bm{\phi} \gets \bm{\phi}^{(i+1)}$ 
    \State $\lambda \gets \lambda^{(i+1)}$ 
   \Else
     \State $i \gets i+1$
   \EndIf\label{outerendif}
\EndWhile\label{outerendwhile}
\State \textbf{return} $\bm{\phi}$, $\lambda$
\EndFunction
\end{algorithmic}
\end{algorithm}

%-----------------------------------------------------------------
%-----------------------------------------------------------------
\section{Projection-based Reduced-Order Modelling}\label{ROM} 

This section briefly describes the key stages of constructing a projection-based reduced-order model, including explanations of proper orthogonal decomposition (POD) and the method of snapshots, singular value decomposition, selecting POD basis functions and constructing a projection-based reduced-order model.

\subsection{Proper Orthogonal Decomposition}\label{ROM-POD}
Proper orthogonal decomposition  was introduced by Lumley~\cite{Lumley1967} to identify the most energetic coherent structures within turbulent flow fields. The method can be applied to data from experiments as well as from numerical simulations. In the latter case, the numerical data is gathered into a matrix, $\bm{S}$, whose columns are individual snapshots pertaining to the solution of the high-fidelity model for a particular set of parameters. For a problem with $N$ spatial degrees of freedom and $M$ snapshots, the snapshots' matrix will be of size $N$ by $M$. According to POD, a snapshot can be decomposed into a finite number of coefficients ($\alpha_i$) and POD basis functions or modes ($\bm{\psi}_i$) as follows:
\begin{equation}\label{eq:POD}
    \bm{\phi} = \sum_{i}\alpha_i \bm{\psi}_i\,,
\end{equation}
where $\alpha_i$ is a scalar and $\bm{\phi},\,\bm{\psi}_i\in\mathbb{R}^N$. POD gives the optimal linear representation of the data in the snapshots' matrix. In order to find the POD basis functions, the average squared error between the snapshots and their projection onto the basis functions is minimised~\cite{holmes2012turbulence}. This leads to an eigenvalue problem, which was originally written in the form 
\begin{equation}\label{SST}
    \bm{S}\bm{S}^T \bm{\psi}_j = \mu_j \bm{\psi}_j\quad\forall j\in N\,,
\end{equation}
where $\mu_j$ is the eigenvalue, $\bm{\psi}_j$ is the eigenvector, and, without loss of generality, it is assumed that the eigenvalues are given in descending order and that the eigenvectors have been orthonormalised. For a large number of degrees of freedom, this problem quickly becomes challenging to solve, and it was Sirovich~\cite{sirovich1987turbulence} who noted that the same (non-zero) eigenvalues could be found by instead solving the following, more tractable problem
\begin{equation}\label{STS}
    \bm{S}^T\bm{S} \bm{\varphi}_k = \mu_k \bm{\varphi}_k\quad\forall k\in M\,,
\end{equation}
for eigenvectors $\bm{\varphi}_k$. Called the method of snapshots, its lower computational burden has led to a wide uptake of POD for analysing results from numerical simulations~\cite{Taira2017,brunton2020machine}. Assuming, again, that the eigenvalues are arranged in descending order, the eigenvectors of the two problems~\eqref{SST} and~\eqref{STS} are related through
\begin{equation}\label{STS_POD_functions}
    \bm{\psi}_k = \frac{1}{\sqrt{\mu_k}}\bm{S}\bm{\varphi}_k \quad \forall k \in M\,. 
\end{equation}
Finding the POD basis functions can also be done by taking a singular value decomposition of the snapshots' matrix~\cite{holmes2012turbulence}, which we now describe briefly. For a real \(N \times M\) matrix $S$, where \(N \geqslant M\), the singular value decomposition of $\bm{S}$ is as follows: 
\begin{equation}\label{eq:SVD}
\bm{S} = \bm{U} \bm{\Sigma} \bm{V^T},
\end{equation}
where $\bm{U}$ is an $N \times N$ matrix consisting of orthonormalised eigenvectors associated with the $M$ largest eigenvalues of $\bm{SS}^T$ and $\bm{V}$ is an $M \times M$ matrix consisting of orthonormalised eigenvectors associated with $\bm{S}^T\bm{S}$. Here, the POD basis functions are the columns of $\bm{U}$. The matrix $\Sigma$ contains the non-negative square roots of the eigenvalues of $\bm{S}^T\bm{S}$ (called the singular values) on its diagonal. These singular values are ordered as follows: 
\begin{equation}
\sigma_1 \geqslant \sigma_2 \geqslant \cdots \geqslant \sigma_M \geqslant 0\,.
\end{equation}
A reduced rank approximation to $\bm{S}$ can be found by setting to zero all but the first~$P$ largest singular values in the matrix~$\Sigma$, and then pre- and post-multiplying by $\bm{U}$ and $\bm{V}^T$ respectively. This will give the optimal rank~$P$ approximation of the matrix~$\bm{S}$ according to the Eckart-Young-Mirsky theorem~\cite{Eckart1936}. 

\subsection{POD Basis Functions}
For some problems, a few basis functions will capture much of the behaviour of the system (as represented by the snapshots), and this is seen by the magnitude of the first few singular values being much larger than the remaining singular values. In such cases, the POD basis can be truncated without introducing much error. To determine how many basis functions should be kept, an empirical expression can be used that relates the number of basis functions retained to the fraction of information captured by them. The amount of information carried by each basis function depends on the square of its singular value. Suppose that the fraction of information of the original system to be captured is $\gamma$, where $0 \leqslant \gamma \leqslant 1$, then the lowest integer value of~$P$ is sought, such that the following is satisfied:
\begin{equation}
\frac{\sum^P_{i=1}\sigma_i^2}{\sum^M_{i=1}\sigma_i^2} \geqslant \gamma\,.
\end{equation}
After truncation is applied (if indeed it is), the basis functions are stored in the matrix~$\bm{R}\in\mathbb{R}^{N\times P}$. If the SVD has been used then the first~$P$ basis functions from $\bm{U}$ are retained to make up the columns of the matrix $\bm{R}$. If no truncation is applied then $P=M$ (i.e.~the number of basis functions is equal to the number of snapshots). If either of the related eigenvalue problems are solved to determine the basis functions, then, once orthonormalised,  the first $P$ eigenvectors are retained to make up $\bm{R}$, either from Equation~\eqref{SST} or Equations~\eqref{STS} and~\eqref{STS_POD_functions}. The complexity of solving the eigenvalue problems is $\mathcal{O}(N^3)$ for Equation~\eqref{SST} and $\mathcal{O}(M^3)$ for Equation~\eqref{STS}. For an SVD, the complexity is $\mathcal{O}(NM^2)$ (see Golub et al.~\cite{Golub1996}), thus, for large problems (many snapshots but where $N\gg M$), the eigenvalue problem in Equation~\eqref{STS} will be the cheapest to solve.

Once $\bm{R}$ has been determined, the $P$ reduced variables $\bm{\alpha}$ associated with a snapshot~$\bm{\phi}$ can be determined from the basis functions by: 
\begin{equation}\label{alpha=rt}
    \bm{\alpha}=\bm{R^T}\bm{\phi},
\end{equation}
and the cell-based scalar flux values can be recovered from the reduced variables by:
\begin{equation}\label{sigma=Ralpha}
    \bm{\phi}=\bm{R}\bm{\alpha}\,,
\end{equation}
where $\bm{R}$ represents the POD basis functions and $\bm{\alpha}$ the POD coefficients or reduced variables. Equation~\eqref{sigma=Ralpha} is the matrix equivalent of Equation~\eqref{eq:POD}.

\subsection{Constructing a Reduced-Order Model for Criticality with POD}\label{sec:PODGROM}
Having found the POD basis functions for the system, the discretised governing equations can be projected onto the reduced space. By inserting Equation~\eqref{sigma=Ralpha} into Equation~\eqref{disc_fiss} and pre-multiplying by $\bm{R}^T$, the reduced system of equations is found: 
\begin{equation}\label{pod_eig}
    \bm{R^T}\bm{A}\bm{R}\bm{\alpha}=\bm{\lambda} \bm{R^T} \bm{B} \bm{R}\bm{\alpha}\,,
\end{equation}
%equation~\eqref{disc_fiss}. 
where the matrices $\bm{A}$ and $\bm{B}$ both depend on the material parameters $\Sigma^a$, $\Sigma^s$, $\Sigma^f$. This reduced system of equations has matrices of size~$P$ by~$P$, whereas the original system given in Equation~\eqref{disc_fiss} is~$N$ by~$N$ where $N\gg P$. 

In this paper, the aim is to construct a reduced-order model based on a set of snapshots, each of which correspond to a particular control-rod configuration, and to use this model to predict solutions for previously unseen control-rod configurations. %
An important point to address is how to modify the reduced system in Equation~\eqref{pod_eig} for unseen parameters. One possibility is to exploit an affine decomposition of the matrices~$\bm{A}$ and~$\bm{B}$, 
by calculating $\bm{R}^T\bm{A}_i\bm{R}$ and $\bm{R}^T\bm{B}_i\bm{R}$ as part of the offline stage, where $\bm{A}_i$ is the $i$th matrix in the affine decomposition of $\bm{A}$, similarly for $\bm{B}_i$. These matrices could be interpolated in order to approximate the matrices corresponding to an unseen parameter~\cite{Heaney2018}. %\textcolor{red}{Although accurate for the problem solved in Heaney et~al.~\cite{Heaney2018}, this method can in general be inaccurate~\cite{Constantine}}. 
Another approach would be to use the Matrix Discrete Empirical Interpolation Method~\cite{Wirtz2014}, which involves sampling the high-fidelity model matrices~$\bm{A}$ and~$\bm{B}$ at a relatively low number of points  before pre- and post-multiplying by $\bm{R}^T$ and $\bm{R}$ respectively, in order to estimate~$\bm{A}$ and~$\bm{B}$ for the unseen parameter. For both these methods, the sampling of $\bm{A}$ and $\bm{B}$, and their pre- and post-multiplication by the POD basis functions would be part of the offline stage, whereas the interpolation of the resulting reduced matrices would be part of the online stage. As the aim of this paper is to demonstrate an autoencoder-based ROM for eigenvalue problems, a simple method is chosen which avoids approximation due to sampling, but would be impractical for large problems. To evaluate $\bm{R}^T\bm{AR}$ and $\bm{R}^T\bm{BR}$ for an unseen parameter, the high-fidelity model is used to assemble~$\bm{A}$ and~$\bm{B}$ (Equations~\eqref{eq:discretised_hfm} and~\eqref{disc_fiss}) and then pre- and post-multiplied by $\bm{R}^T$ and $\bm{R}$ respectively. This would not be practical for systems with large number of degrees of freedom, as the online stage, in which $\bm{R}^T\bm{AR}$ and $\bm{R}^T\bm{B}\bm{R}$ are approximated for the unseen parameter, should be independent of the high-fidelity model for reasons of computational efficiency. Future work will involve larger scale problems, for which more computationally efficient methods will be investigated.
%Another alternative would be to use the Discrete Empirical Interpolation method~\cite{Chaturantabut2010,Antil2014}. 

For the POD-based ROM, the offline stage consists of solving the high-fidelity model many times for different material parameters (cross-sections), see Equation~\eqref{disc_fiss} and Algorithm~\ref{alg:PM_outer}. In this paper, the POD basis functions are determined by solving the eigenvalue problem associated with the method of snapshots, see Equation~\eqref{STS}. The online stage of the POD-based ROM consists of assembling the matrices~$\bm{A}$ and~$\bm{B}$ for a given set of parameters, projecting these matrices onto the reduced space and then solving with the power method. The outer iteration algorithm for the POD-based ROM is very similar to that of the high-fidelity model (Algorithm~\ref{alg:PM_outer}), except for the passing down of the POD basis functions from the outer to the inner iterations. The so-called inner iterations are described for POD-based ROM in Algorithm~\ref{alg:PODROM_inner}. (Although there are no iterations, the description `inner iterations' is still used for Algorithm~\ref{alg:PODROM_inner}). Here the reduced-order model, Equation~\eqref{pod_eig}, is used to solve for the scalar fluxes. The right-hand side of Equation~\eqref{pod_eig} is evaluated using the current value of scalar flux, resulting in a system that can be solved for the reduced variables. As the system is small, the problem is solved directly with Gaussian elimination. There is no need to use the eigenvalue in the source term (see line~\ref{alg2_source_line}) as, once passed back to the outer iterations, the flux is normalised. However, it is included here as it is needed in the inner iterations for the autoencoder-based reduced order models.

\begin{algorithm}
\caption{POD-based reduced-order model: inner iterations}\label{alg:PODROM_inner}
\begin{algorithmic}[1]
\Function{POD\textunderscore{}Inner\textunderscore{}Iterations}{$\bm{A}$, $\bm{B}$, $\bm{\phi}$, $\lambda$, $\bm{R}$} 
\State $\bm{s} = \lambda \left( \bm{R}^T \bm{B} \right) \bm{\phi}$ \label{alg2_source_line}
\Comment{set a source with values from the outer iterations}
\State solve the reduced-order model for $\bm{\alpha}$:
\State $\left( \bm{R}^T \bm{AR}\right)\bm{\alpha} = \bm{s} $ %\lambda \left( \bm{R}^T \bm{B} \bm{R} \right) \bm{\phi}$
\State $\bm{\phi} = \bm{R\alpha}$ \Comment{find the updated scalar flux from the reduced variables}
\State \textbf{return} $\bm{\phi}$
\EndFunction
\end{algorithmic}
\end{algorithm}

%-----------------------------------------------------------------
%-----------------------------------------------------------------
\section{Autoencoders and their inclusion in projection-based ROMs}\label{AEs} 
First, the two types of autoencoder used in this paper are described: a standard, fully-connected autoencoder and an SVD-autoencoder. Then, their incorporation into projection-based reduced-order models is explained. The use of autoencoders for dimensionality reduction in eigenvalue problems represents the main novelty in this paper. To the authors' best knowledge, the hybrid SVD-autoencoder is also a new method.

\subsection{Autoencoders}
\subsubsection{A Standard Autoencoder}\label{AE}
An autoencoder is a special type of feed-forward neural network that is trained to learn the identity map. A bottleneck at its central layer forces the autoencoder to learn a reduced representation of the data. The simplest (vanilla) autoencoder has an input layer with $n$~neurons, a hidden layer with $m$ neurons where $m<n$, and an output layer with $n$~neurons. The values of the neurons in the hidden layer define the latent space and can be referred to as the latent variables. Networks with additional hidden layers can learn more complicated patterns creating a more effective representation of complex data. An autoencoder can be split into two networks: the encoder, which maps an input to the latent variables, and a decoder, which maps from the latent variables to the output. If the input vector is given by $\bm{x}\in\mathbb{R}^n$, the output vector by $\hat{\bm{x}}\in\mathbb{R}^n$ and the latent variables of the hidden layer by $\bm{y}\in\mathbb{R}^m$, then the autoencoder can be written as a composition of functions as follows
\begin{equation}\label{eq:autoencoder_mapping}
\hat{\bm{x}} = f^{\text{AE}}(\bm{x}; \bm{w}^{\text{AE}}) = f^{\text{DEC}}\left( \bm{y; \bm{w}^{\text{DEC}}} \right) = f^{\text{DEC}}\left(  f^{\text{ENC}}(\bm{x}; \bm{w}^{\text{ENC}}) ; \bm{w}^{\text{DEC}}\right) 
\end{equation}
where $f^{\text{AE}}$, $f^{\text{ENC}}$ and $f^{\text{DEC}}$ represent the autoencoder, encoder and decoder respectively, with associated weights $\bm{w}^{\text{AE}}$, $\bm{w}^{\text{ENC}}$ and $\bm{w}^{\text{DEC}}$.

\begin{figure}[!htb]
\def\layersep{2.5cm}

\begin{tikzpicture}[shorten >=1pt,->,draw=black!50, node distance=\layersep]
    \tikzstyle{every pin edge}=[<-,shorten <=1pt]
    \tikzstyle{neuron}=[circle,fill=black!25,minimum size=17pt,inner sep=0pt]
    \tikzstyle{input neuron}=[neuron, fill=darkgreen];
    \tikzstyle{output neuron}=[neuron, fill=darkred];
    \tikzstyle{hidden neuron1}=[neuron, fill=darkblue];
    \tikzstyle{hidden neuron2}=[neuron, fill=darkblue];
    \tikzstyle{annot} = [text width=4em, text centered]

    % Draw the input layer nodes
    \foreach \name / \y in {1,...,6}
    % This is the same as writing \foreach \name / \y in {1/1,2/2,3/3,4/4}
        \node[input neuron] (I-\name) at (0,-\y) {};

    % Draw the hidden layer nodes
    \foreach \name / \y in {1,...,4}
        \path[yshift=-1cm]
            node[hidden neuron1] (H1-\name) at (\layersep,-\y cm) {};
% Draw the hidden layer nodes
    \foreach \name / \y in {1,...,2}
        \path[yshift=-2cm]
            node[hidden neuron2] (H2-\name) at (2*\layersep,-\y cm) {};
    % Draw the output layer node
    %\node[output neuron,pin={[pin edge={->}]right:Output}, right of=H-3] (O) {};
    \foreach \name / \y in {1}
    \path[yshift=-2.5cm]
    node[output neuron] (O) at (3*\layersep,-\y cm) {};

    % Connect every node in the input layer with every node in the
    % hidden layer.
    \foreach \source in {1,...,6}
        \foreach \dest in {1,...,4}
            \path (I-\source) edge (H1-\dest);
    
    \foreach \source in {1,...,4}
        \foreach \dest in {1,...,2}
            \path (H1-\source) edge (H2-\dest);

    % Connect every node in the hidden layer with the output layer
    \foreach \source in {1,...,2}
        \path (H2-\source) edge (O);

    % Draw the hidden layer nodes
    \foreach \name / \y in {1,...,2}
        \path[yshift=-2cm]
            node[hidden neuron1] (H3-\name) at (4*\layersep,-\y cm) {};
\foreach \name / \y in {1,...,4}
        \path[yshift=-1cm]
            node[hidden neuron2] (H4-\name) at (5*\layersep,-\y cm) {};
    % Draw the output layer nodes
    \foreach \name / \y in {1,...,6}
    % This is the same as writing \foreach \name / \y in {1/1,2/2,3/3,4/4}
        \node[input neuron] (I2-\name) at (6*\layersep,-\y) {};

    % Connect every node in the hidden layer with the output layer
    \foreach \dest in {1,...,2}
        \path (O) edge (H3-\dest);

    % Connect every node in the input layer with every node in the
    \foreach \source in {1,...,4}
        \foreach \dest in {1,...,2}
            \path (H3-\dest) edge (H4-\source);
    \foreach \source in {1,...,6}
        \foreach \dest in {1,...,4}
            \path (H4-\dest) edge (I2-\source);

    % Annotate the layers
%    \node[annot,above of=H-1, node distance=1cm] (hl) {Hidden layer};
%    \node[annot,left of=hl] {Input layer};
%    \node[annot,right of=hl] {Output layer};

    \node[annot,above of=I-1, node distance=1cm] (input) {Input};
    \node[annot,right of=input] (encoder1) {Encoder};
    \node[annot,right of=encoder1] (encoder2) {Encoder};
    \node[annot,right of=encoder2] (latent) {Latent Variable(s)};
    \node[annot,right of=latent] (decoder1) {Decoder};
    \node[annot,right of=decoder1] (decoder2) {Decoder};
    \node[annot,right of=decoder2] {Output};

\end{tikzpicture}
% End of code
\caption{The architecture of a fully-connected autoencoder comprising of three encoder layers and three decoder layers. This network compresses 6 variables into 1 latent variable and then expands back out to the original 6 variables.}
\label{fig:autoencoder}
\end{figure}
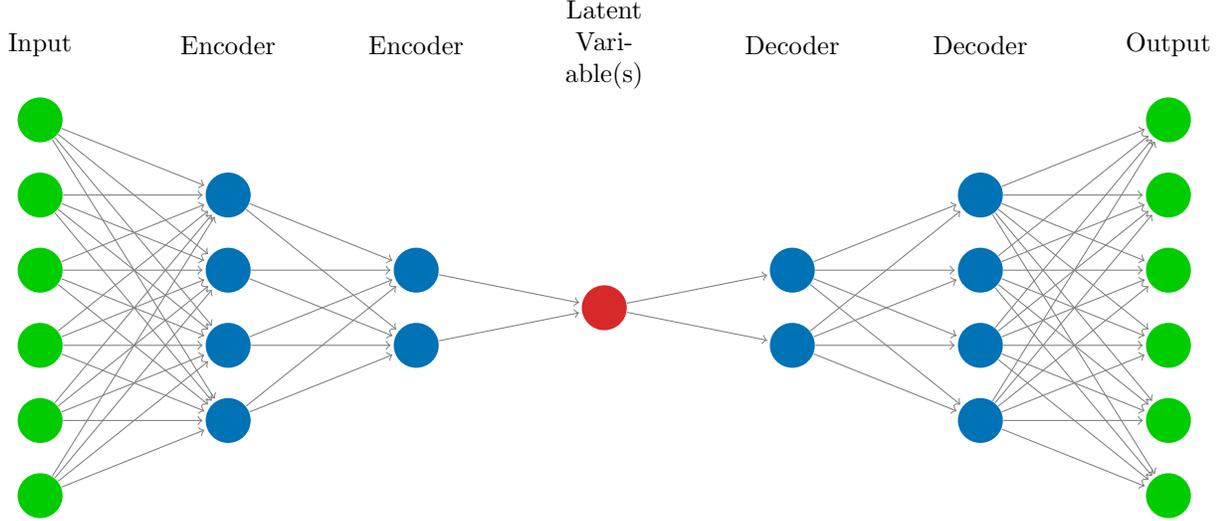

Figure~\ref{fig:autoencoder} illustrates a fully-connected autoencoder with 6~layers (only layers with trainable weights are included in the total), 5 of these are hidden layers. Neurons can have connections with sending neurons from the previous layer (except for those in the input layer) and with receiving neurons from the following layer (except for those in the output layer). In this diagram, the input and output layers have 6~neurons each and the latent space is of dimension~1. The grey arrows indicate how the information passes through the autoencoder from left to right. The input to a neuron is calculated as a weighted sum of the outputs of the sending neurons to which it is connected and a bias term. The output of the neuron is the evaluation of the activation function for this weighted sum. For example, the output of the $i$th neuron of the hidden layer of the vanilla autoencoder described above, is given by
\begin{equation}
y_i = g\left( \sum_{j=1}^n w^{\text{ENC}}_{ij} x_j + b_i \right)
\end{equation}
where $g$ is the activation function, $b_i$ is the bias of the $i$th neuron, $\{x_j\}_{j=1}^n$ are the values of the sending neurons connected to the $i$th neuron and $w_{ij}$ is the weight of the $j$th sending neuron of the $i$th neuron in the hidden layer. 

The data available to build the network is customarily split into three parts: training data, validation data and test data. The use of these data sets in the training process will now be described. The activation function and number(s) of neurons in the hidden layer(s) are examples of hyper-parameters of the network, and the weights and biases are referred to as parameters of the network. For given hyper-parameters, the objective of training the autoencoder is to find weights that minimise the so-called loss function, which, for regression problems, is often taken as the mean squared error of the reconstruction:
\begin{equation}
\frac{1}{N^{\text{ex}}}\sum_{j=1}^{N^{\text{ex}}} || \bm{x}_j - \hat{\bm{x}}_j||^2 = \frac{1}{N^{\text{ex}}}\sum_{j=1}^{N^{\text{ex}}} || \bm{x}_j - f^{\text{ae}}\left(\bm{x}_j; \bm{w}^{\text{ae}}\right)||^2\ ,
\end{equation}
where $N^{\text{ex}}$ is the number of samples or examples used to train the autoencoder and the norm used is the Euclidean norm. Based on gradient descent-type methods, a back-propagation algorithm is used to optimise the weights by solving
\begin{equation}
\bm{w}^{\text{ae}} = \argmin_{\tilde{\bm{w}}^{\text{ae}}} \frac{1}{2}\sum_{j=1}^{N^{\text{ex}}} || \bm{x}_j - f^{\text{ae}}\left(\bm{x}_j; \tilde{\bm{w}}^{\text{ae}}\right)||^2\ .
\end{equation}
In order to find the optimal hyper-parameters, less sophisticated approaches than back-propagation are employed, such as cross-validation~\cite{Kuhn}. This involves training a number of networks each with a different set of hyper-parameters and comparing the performance of these networks on the validation data. The best-performing network is selected and can be re-trained with both the training and validation data combined. The performance of this network is then assessed on the test data.

\subsubsection{SVD Autoencoder}\label{SVD_AE} 
The length of the input to the autoencoder is one factor that governs the difficulty of training an autoencoder.  The longer these vectors are, the more parameters (weights) there are to be adjusted during optimisation, which increases the computational cost. To combat this, the dimensionality reduction can be split into two steps. The first step is to apply an SVD to the inputs of the autoencoder (i.e.~the snapshots). A basis is obtained and the snapshots are projected onto the resulting low-dimensional space yielding reduced variables~$\bm{\alpha}$ (just as is done in POD-based ROM), see Equation~\eqref{alpha=rt}. This reduces the length of the input vectors from the number of degrees of freedom of the problem ($N$) to the number of basis functions ($P$ where $P\ll N$). The second step is to train a fully-connected autoencoder with the reduced variables of the snapshots (or the subset of these that have been selected for training), thereby further reducing the dimensionality.

%\label{Proj_developed} 
\subsection{Projection-based ROMs with autoencoders}\label{Proj_intr} 
POD defines a linear map between the scalar fluxes $\bm{\phi}$, and the corresponding reduced variables, $\bm{\alpha}$, see Equation~\eqref{alpha=rt}. Using an autoencoder for the dimensionality reduction means that the map between these variables will, in general, be nonlinear. Consequently, even for a linear system of equations, the resulting ROM will be nonlinear, and for an eigenvalue problem, this adds an additional level of non-linearity. To construct a reduced-order model, the nonlinear mapping between the reduced variables and the scalar fluxes needs to be linearised and the reduced system of equations needs to be derived. Both points are now elucidated. 

Suppose that we have a known state of reduced variables $\tilde{\bm{\alpha}}$ and scalar fluxes $\tilde{\bm{\phi}}$, and wish to find the change in the scalar flux $\Delta\bm{\phi}$ due to a small change in the reduced variables $\Delta\bm{\alpha}$. We can linearise the map $\bm{C}$ between these variables using a first order Taylor series as follows:
\begin{eqnarray}
\Delta\bm{\phi}                   & = & \bm{C}(\tilde{\bm{\alpha}})\, \Delta\bm{\alpha} \,, \label{eq:C_linearised} \\
\textrm{where}\ \Delta\bm{\phi}   & = & \bm{\phi} - \tilde{\bm{\phi}}\,, \\
              \ \Delta\bm{\alpha} & = & \bm{\alpha} - \tilde{\bm{\alpha}} \,, \\
\textrm{and}\ \bm{C}(\tilde{\bm{\alpha}})   & = &  \frac{\textrm{d}\bm{\phi}}{\textrm{d}\bm{\alpha}} \bigg|_{\tilde{\bm{\alpha}}} \label{eq:C} \,. 
\end{eqnarray}
The scalar flux $\bm{\phi}$ can therefore be approximated as
\begin{equation}\label{eq:phi_approximation}
\bm{\phi} = \bm{C}(\tilde{\bm{\alpha}})\bm{\alpha} - \bm{C}(\tilde{\bm{\alpha}})\tilde{\bm{\alpha}} + \tilde{\bm{\phi}}\ .
\end{equation}

In order to calculate $\bm{C}(\tilde{\bm{\alpha}})$, perturbations are made to $\tilde{\bm{\alpha}}$ from which  the entries of $\bm{C}$ can be inferred. Each of the $P$ entries in $\tilde{\bm{\alpha}}$ are perturbed in turn by a small value, which yields a column of entries in $\bm{C}$. The $k$th perturbation simply adds a small number~$\epsilon$ to the $k$th component of $\tilde{\bm{\alpha}}$, so, for example, the second perturbation is given by
\begin{equation}
{}^{2}\Delta\bm{\alpha} =  \underbrace{ \left(  0, \epsilon, 0, 0, \ldots, 0 \right)}_{P\ \text{entries}}{}^T \ .
\end{equation}

The decoder can be used to calculate the scalar flux corresponding to the perturbed reduced variables $\bm{{}^k\alpha}$:
\begin{equation}
{}^k\bm{\phi} = f^{\text{DEC}}\left( {}^k\bm{\alpha} \right) = f^{\text{DEC}}\left( \tilde{\bm{\alpha}} + {}^k\Delta\bm{\alpha}\right).
\end{equation}
The change in the scalar flux can be calculated by subtracting the known state $\tilde{\bm{\phi}}$ from the above. The $k$th column of $\bm{C}$ is now known:
\begin{equation}\label{eq:calculate_C}
C_{jk}(\tilde{\bm{\alpha}}) =
\frac{
{}^k\phi_j - \tilde{\phi}_j 
}{\epsilon} 
=  \frac{{}^k\Delta\phi_j}{\epsilon} \qquad \forall j\ ,
\end{equation}
where $\phi_j$ is the $j$th component of the scalar flux. Once every component of the reduced variables $\bm{\alpha}$ has been perturbed, the $\bm{C}$ matrix will be fully determined. Having found the map between the reduced variables and the high-fidelity model variables, the projection-based ROM, using an autoencoder, can now be developed. 

When solving the high-fidelity model, an inner iterative loop solves for the scalar flux and an outer iterative loop solves for the eigenvalue. When using a reduced-order model with an autoencoder, we solve a reduced-order model for the scalar flux but use the same outer iterations to solve for the eigenvalues, similar to what was done for the POD-based method. Hence, Algorithm~\ref{alg:PM_outer} is used unchanged for the reduced model. The right-hand side of Equation~\eqref{disc_fiss} is treated as a source which will be held constant while solving, iteratively, for the flux. In other words, although $\bm{\phi}$ and $\lambda$ will be updated in the inner iterations, the source term will use the values of these variables from the outer iterations:
\begin{equation}\label{eq:full_gevp_source_rhs}
\bm{A\phi} = \bm{s} \,,
\end{equation}
where $\bm{s} = \lambda^{(\cdot)} \bm{B} \bm{\phi}^{(\cdot)}$.
Substituting Equation~\eqref{eq:phi_approximation} into Equation~\eqref{eq:full_gevp_source_rhs}, rearranging and pre-multiplying the result by $\bm{C}^T$ gives the reduced system of equations for the autoencoder-based ROM:
\begin{equation}
\left( \bm{C}^T(\tilde{\bm{\alpha}}) \bm{A}  \bm{C}(\tilde{\bm{\alpha}}) \right) \Delta\bm{\alpha} = \bm{C}^T(\tilde{\bm{\alpha}}) \bm{s} - \bm{C}^T(\tilde{\bm{\alpha}}) \bm{A} \tilde{\bm{\phi}}\,.
\end{equation} 

As this system is relatively small for the examples used in this paper, it is solved with Gaussian elimination. The algorithm for the inner iterations, which solve for the scalar flux, is given in Algorithm~\ref{alg:ROM_inner} in which the dependence of $\bm{C}$ on $\tilde{\bm{\alpha}}$ has been dropped for readability. Regularisation has been applied to the reduced-order model by adding a small number, $\varepsilon$, to the diagonal of $\bm{C}^T\bm{AC}$ as, without this, there is no guarantee that this matrix will be non-singular. The regularisation will also ensure that the iterative change in reduced variables ($\Delta\bm{\alpha}$) will not become too large. Note that this is a modified power iteration, as the eigenvalue is included in the source term. For the standard power method this is not required when determining the flux, as normalisation renders this unnecessary. However, in line~\ref{alg:line9alg3} of the algorithm, the right-hand side of this expression has a source term and an additional term, meaning that the size of the source term is now important. Both autoencoders described in Sections~\ref{AE} and~\ref{SVD_AE} can be used in the reduced-order model framework described in this section.

\begin{algorithm}
\caption{AE-based reduced-order model: inner iterations}\label{alg:ROM_inner}
\begin{algorithmic}[1]
\Function{Inner\textunderscore{}Iterations}{$\bm{A}$, $\bm{B}$, $\bm{\phi}$, $\lambda$} 
\State $k$\textunderscore{}$max = 100$
\State $\alpha$\textunderscore{}$tol = \num{e-8}$
\State $\bm{s} = \lambda \bm{B} \bm{\phi}$
\Comment{set a source with values from the outer iterations}
\State $\bm{\phi}^{(0)} = \bm{\phi}$; $\lambda^{(0)} = \lambda$; $\bm{\alpha}^{(0)} = f^{\text{ENC}}(\bm{\phi}^{(0)})$ 
\Comment{initialise with values from the outer iterations}
\State $k = 0$
\State \texttt{not\textunderscore{}converged = True}
\While{\texttt{not\textunderscore{}converged}}
  \State calculate $\bm{C}^{(k)}$ for $\bm{\alpha}^{(k)}$ using Equation~\eqref{eq:calculate_C}
  \State solve the reduced-order model:
  \State $\left( (\bm{C}^{(k)})^T \bm{A} \bm{C}^{(k)} + \varepsilon\bm{I} \right)\Delta\bm{\alpha}^{(k+1)} =  (\bm{C}^{(k)})^T \bm{s} - (\bm{C}^{(k)})^T \bm{A} \bm{\phi}^{(k)}$ \label{alg:line9alg3} %\Comment{the reduced-order model}
  \State $\bm{\alpha}^{(k+1)} = \bm{\alpha}^{(k)} + \Delta\bm{\alpha}^{(k+1)}$ \Comment{update variables}
  \State $\Delta\bm{\phi}^{(k+1)} = \bm{C}^{(k)}\Delta\bm{\alpha}^{(k+1)}$
  \State $\bm{\phi}^{(k+1)} = \bm{\phi}^{(k)} + \Delta\bm{\phi}^{(k+1)}$ 
  \If{ ($\Delta\bm{\alpha}^{(k+1)}<\alpha$\textunderscore{}$tol$ \textbf{or} $k=k$\textunderscore{}$max$)}
    \State \texttt{not\textunderscore{}converged = False}
    \State $\bm{\phi} \gets \bm{\phi}^{(k+1)}$  
  \Else
    \State $k \gets k+1$
  \EndIf\label{euclidendif}
\EndWhile\label{innerendwhile}
\State \textbf{return} $\bm{\phi}$
\EndFunction
\end{algorithmic}
\end{algorithm}

\section{Results}\label{results_disc} 

The methods outlined in Sections~\ref{ROM} (POD-based ROM) and~\ref{AEs} (autoencoder-based ROMs) are applied to two test cases. Before presenting the results, the error measures that are used to compare the methods are introduced, followed by an explanation of how the material parameters (cross-sections) are homogenised. The first test case is a 1D slab reactor and results from a POD-based ROM and an AE-based ROM % and a VAE-based ROM 
are presented. Further comparisons are made between a POD-based ROM and an AE-based ROM constructed with fewer reduced variables. The second test case is a simplified 2D reactor. For this case, results are presented for ROMs based on POD, on an autoencoder and a hybrid SVD-autoencoder. The reduced-order models were developed in python using Keras~\cite{chollet2015keras} with the TensorFlow backend.

\subsection{Error measures}\label{sec:error_definitions}
To assess the accuracy of the results, a number of error measures are introduced. In all the expressions, the reduced-order model solutions are measured against the high-fidelity model results. First, the error committed in the dimensionality reduction step is evaluated. For the POD-based ROM, this error is incurred by projecting the high-fidelity model solutions onto the reduced space and for the AE-based ROM this error occurs during the compression and de-compression of the high-fidelity model solutions by the autoencoder. For the SVD-AE-based ROM, the error is due to both projecting the solution onto a space defined by POD basis functions, and compressing and de-compressing with an autoencoder. The normalised maximum errors in the reconstruction of the solution are defined here  using
\begin{align} 
    e_{\text{max}}^{\text{POD}}(\bm{\phi}^{\text{HFM}})&=\frac{ \phi^{\text{HFM}}_k-\left(\bm{R}\bm{R}^T\bm{\phi}^{\text{HFM}} \right)_k}{\| \bm{\phi}^{\text{HFM}}\|_\infty}\, \quad &
    {\text{where}}\;  k & = \argmax_{i\,\in\, \text{cells}} \mid\phi^{\text{HFM}}_i-\left(\bm{R}\bm{R}^T\bm{\phi}^{\text{HFM}} \right)_i\mid\ ,
    \label{error_compression_POD}\\[2mm]
    e_{\text{max}}^{\text{AE}}(\bm{\phi}^{\text{HFM}})&=\frac{ \phi^{\text{HFM}}_k-\left(f^{\text{AE}}(\bm{\phi}^{\text{HFM}})\right)_k}{\| \bm{\phi}^{\text{HFM}}\|_\infty}\, \quad &
    {\text{where}}\;  k &= \argmax_{i\,\in\, \text{cells}} \mid\phi^{\text{HFM}}_i-\left(f^{\text{AE}}(\bm{\phi}^{\text{HFM}})\right)_i\mid\ ,\label{error_compression_AE}\\[2mm]
    e_{\text{max}}^{\text{SVD-AE}}(\bm{\phi}^{\text{HFM}}) & = \frac{ \phi^{\text{HFM}}_k-\left(\bm{R} f^{\text{AE}}\left(\bm{R}^T\bm{\phi}^{\text{HFM}}\right) \right)_k}{\| \bm{\phi}^{\text{HFM}}\|_\infty} \quad & \text{where} \;  k & = \argmax_{i\,\in\, \text{cells}} \mid\phi^{\text{HFM}}_i-\left(\bm{R} f^{\text{AE}}\left(\bm{R}^T\bm{\phi}^{\text{HFM}}\right))\right)_i\mid\ , \label{error_reconstruction_SVDAE}
\end{align}
for the POD-based ROM, the AE-based ROM and SVD-AE-based ROM respectively, where $\|\cdot\|_\infty$ represents the maximum norm, $i$ and $k$ represent cell indices, and the superscript `$\text{HFM}$' indicates the flux is associated with the high-fidelity model.

Second, the error in the scalar flux of the reduced-order models is calculated. The normalised maximum error in the flux profile is determined by
\begin{equation}\label{error}
    e_{\text{max}}(\bm{\phi}^{\text{ROM}})=\frac{ \phi^{\text{HFM}}_j-\phi_j^{\text{ROM}}}{\| \bm{\phi}^{\text{HFM}}\|_\infty} \quad \text{where} \; j = \argmax_{i\,\in\, \text{cells}} \mid\phi^{\text{HFM}}_i-\phi^{\text{ROM}}_i\mid\ ,
\end{equation}
in which the superscript `$\text{ROM}$' indicates the solution is from a reduced-order model (either POD-based or autoencoder-based), and $i$ and $j$ represent cell indices.  

Calculating the error in $k_{\text{eff}}$ is done by comparing the value calculated by the high-fidelity model with that of the reduced-order model:
\begin{equation}\label{keff_error}
    e(k_{\text{eff}}^{\text{ROM}})=k_{\text{eff}}^{\text{HFM}}-k_{\text{eff}}^{\text{ROM}}.
\end{equation}
As $k_{\text{eff}}$ is usually close to one, no normalisation is applied.

Finally, an average maximum error is defined, where, for $N$ sets of material parameters resulting in $N$ solutions of the high-fidelity model and $N$ solutions of one of the reduced-order models:
\begin{equation}\label{average_maximum_error}
%    \frac{1}{N}\sum_{i=1}^{N} %|e^{\text{flux}}(\mu_i)|
\bar{e}_{\text{max}}(\bm{\phi}^{\text{ROM}}) = \frac{1}{N}\sum_{l=1}^{N} \frac{ \left|\,\phi^{\text{HFM}}_j(\bm{\mu}_l)-\phi_j^{\text{ROM}}(\bm{\mu}_l)\,\right|}{\| \phi^{\text{HFM}}(\bm{\mu}_l)\|_\infty}
\end{equation}
where $\bm{\mu}_l$ refers to the parameters (i.e.~the macroscopic cross-sections) used for the $l$th problem. Until this point, dependence of the solution on the material parameters has not been explicitly indicated (for readability), however here it is required. A similar error measure can also be used to find the average maximum error in $k_{\text{eff}}$:
\begin{equation}\label{average_maximum_error_keff}
\bar{e}(k_{\text{eff}}^{\text{ROM}}) = \frac{1}{N}\sum_{l=1}^{N} \left|\, k_{\text{eff}}^{\text{HFM}}(\bm{\mu}_l)-k_{\text{eff}}^{\text{ROM}}(\bm{\mu}_l)\,\right|\ ,
\end{equation}
and the average maximum error in the reconstruction:
\begin{equation}\label{average_maximum_error_recon}
\bar{e}_{\text{max}}^{\text{\,DR}} (\bm{\phi}^{\text{HFM}}) = \frac{1}{N}\sum_{l=1}^{N} \left|\, 
e_{\text{max}}^{\text{\,DR}}(\bm{\phi}^{\text{HFM}})
\,\right|\ ,
\end{equation}
where DR represents one of the three dimensionality reduction methods used: POD, AE or SVD-AE.
%\begin{equation}\label{average_maximum_error_keff}
%\frac{1}{N}\sum_{k=1}^{N} \frac{ \|k_{\text{eff}}^{\text{HFM}}(\bm{\mu}_k)-k_{\text{eff}}^{\text{ROM}}(\bm{\mu}_k)\|_\infty}{\| k_{\text{eff}}^{\text{HFM}}(\bm{\mu}_k)\|_\infty}
%\end{equation}

\subsection{Homogenising the Material Parameters}\label{sec:homogenisation} In the test cases that follow, there are fuel regions and also control-rod regions, within which the control rods can be fully inserted, completely withdrawn or partially inserted. In the latter case, the cross-sections will be calculated by combining the cross-sections for fuel and control rod. %
Control-rod regions are assumed to have a uniform distribution of the averaged properties of the mixture within the system based on a mixing coefficient for each region. This can be written as follows
\begin{eqnarray}
\Sigma_{\text{hom}}^a = r \Sigma_{\text{cr}}^a + (1-r) \Sigma_{\text{fuel}}^a \label{eq:mixing_absorption},\\
\Sigma_{\text{hom}}^s = r \Sigma_{\text{cr}}^s + (1-r) \Sigma_{\text{fuel}}^s , \label{eq:mixing_scattering}
\end{eqnarray}
where the mixing coefficient, $r$, for a given region, lies between 0 and 1, and the subscripts `$\text{hom}$', `$\text{cr}$' and `$\text{fuel}$' indicate a homogenised cross-section, and cross-sections of the control rod and the fuel respectively. If $r$ were chosen to be proportional to the amount of control rod present in a control-rod region, the high absorption of the control rods would dominate the homogenised cross-section. In order to address this, the reciprocal of the absorption cross-sections are averaged:
\begin{equation}\label{eq:mixing_reciprocal}
    \frac{1}{\Sigma_{\text{hom}}^a} = \frac{z}{\Sigma_{\text{cr}}^a} + \frac{1-z}{\Sigma_{\text{fuel}}^a}\qquad\text{where}\ z\in[0,1]\,. 
\end{equation}

By combining Equations~\eqref{eq:mixing_absorption} and~\eqref{eq:mixing_reciprocal} it can be shown that the mixing coefficient is given by
\begin{equation}\label{eq:mixing_coefficient}
    r = \frac{z\Sigma_{\text{fuel}}^a}{z\Sigma_{\text{fuel}}^a + (1-z)\Sigma_{\text{cr}}^a}\,.
\end{equation}
In this way, the desired amount of control rod in a region can be chosen by setting $z$, converting this into a value for $r$ by using Equation~\eqref{eq:mixing_coefficient}, and then calculating the homogenised absorption and scattering cross-sections from Equations~\eqref{eq:mixing_absorption} and~\eqref{eq:mixing_scattering}.

\subsection{1D Slab Reactor}\label{1D_slab} 

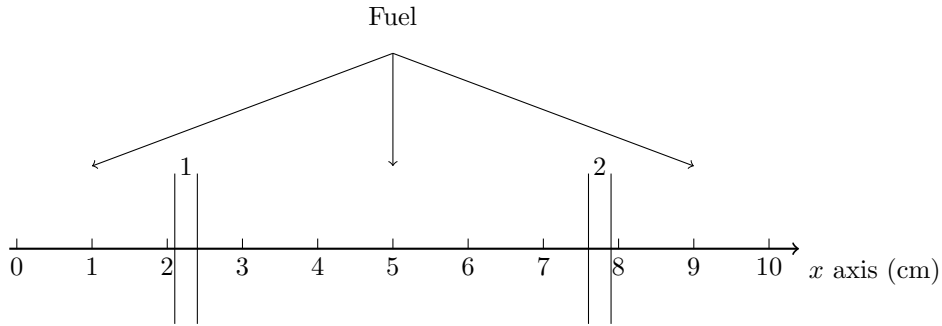
\begin{figure}[!htb]
\def\layersep{2.5cm}
\centering
\begin{tikzpicture}

\draw[thick,->] (3.9,3.9) -- (14.4,3.9) node[anchor=north west] {$x$ axis (cm)};
 \foreach \x/\xtext in {4/0,5/1,6/2,7/3,8/4,9/5,10/6,11/7,12/8,13/9,14/10}
       \draw[yshift=4cm] (\x cm,1pt) -- (\x cm,-3pt) node[anchor=north] {$\xtext$};
\draw (6.1,4.9) -- (6.1,2.9);
\draw (6.4,4.9) -- (6.4,2.9);
\draw (11.9,4.9) -- (11.9,2.9);
\draw (11.6,4.9) -- (11.6,2.9);

\node at (9,7) {Fuel};
\node at (6.25,5) {1};
\node at (11.75,5) {2};
\draw [black, ->          ] (9,6.5) -- (5,5);
\draw [black, ->          ] (9,6.5) -- (9,5);
\draw [black, ->          ] (9,6.5) -- (13,5);
\end{tikzpicture}
\caption{Geometry of a 1D slab reactor consisting  of fuel and two control-rod regions (regions~1 and~2).}
\label{fig:1dcore}
\end{figure}

The length of the slab reactor is \SI{10}{cm}, and it consists of fuel and two control-rod regions labelled as~1 and~2 in Figure~\ref{fig:1dcore}. Region~1 is located between \SI{2.2}{cm} and \SI{2.5}{cm} on the $x$-axis, region~2, between \SI{7.5}{cm} and \SI{7.8}{cm}. This test case was used by Buchan~et~al.~\cite{buchan2013pod}
%long and contains 100 cells and the cells are 1mm wide. 

Table~\ref{1dprop} contains the macroscopic cross-sections used for the control rods and the fuel. The control rods can be fully inserted, completely withdrawn or partially inserted into regions~1 and~2. In the latter case, the cross-sections for these regions will be calculated by combining the cross-sections for fuel and control rod by defining mixing coefficients as described in the previous section. The amount of control rod in region~1 will, in general, be different from that of region~2, henceforth, $z_1$ and $r_1$ refer to region~1, and $z_2$ and $r_2$ refer to region~2. 144~cells were used in the control-volume discretisation although 44~of these were used to enforce the boundary conditions, resulting in a system with 100 degrees of freedom.

\begin{table}[h!]
\centering
 \begin{tabular}{l c c c} 
 \hline\hline
       & $\Sigma^a$ & $\Sigma^s$ & $\Sigma^f$ \\[0.5ex] 
 \hline\hline
 Fuel & 0.45 & 2 & 0.5 \\ 
 %\hline
 Control Rod & 0.9\phantom{5} & 2 & 0\phantom{.5} \\[0.5ex]  
 \hline\hline
\end{tabular}
\caption{Macroscopic cross-sections for the 1D slab reactor (units \si{\per\cm}).}
\label{1dprop}
\end{table}

To generate the data necessary to build the reduced-order models, 200 problems were solved using the high-fidelity model with different control-rod settings. This was done by choosing values for $z_1$ and $z_2$ at random from the interval $[0,1]$ and using Equations~\eqref{eq:mixing_coefficient}, \eqref{eq:mixing_absorption} and~\eqref{eq:mixing_scattering} to calculate the corresponding cross-sections. The control-volume discretisation along with the cross-sections determine the matrices $\bm{A}$ and $\bm{B}$ in Equation~\eqref{disc_fiss}. Given an initial guess for the scalar flux, the system in~\eqref{disc_fiss} is solved with the power method as described in Section~\ref{sec:power_it} yielding 200 solutions of the high-fidelity model. One hundred of the high-fidelity model solutions were treated as snapshots and were used to generate basis functions (either by applying POD or by training an autoencoder). These snapshots are referred to as `seen' data (in the context of ROM) or as training data (in the context of neural networks).  The remaining 100 solutions were used in the online stage to test how good the models were at predicting solutions for unseen data (i.e.~parameter sets whose corresponding solutions were not included in the snapshots). This data is referred to as unseen data (in the context of ROM) and as test data (in the context of  neural networks). Two reduced-order models are compared: a POD-based ROM and an autoencoder-based ROM. % and a VAE-based ROM.
The models have 10~reduced variables (10 POD coefficients for the POD-based model and 10 latent variables for the autoencoder-based model). Following this, a POD-based model and an autoencoder-based model are developed with just two reduced variables.

\subsubsection{POD-based Reduced-Order Model} \label{1D_pod_gal} 
The first method to be tested is a POD Galerkin method which uses POD and the SVD to determine basis functions, projects the discretised governing equations onto the subspace defined by the POD basis functions, and results in a generalised eigenvalue problem of reduced dimension, Equation~\eqref{pod_eig}, to be solved in the online stage for a given set of macroscopic cross-sections. This method is described in Section~\ref{sec:PODGROM}.

\begin{figure}[!htb]
\centering
\begin{minipage}{.45\textwidth}
  \centering
  \includegraphics[scale=0.5]{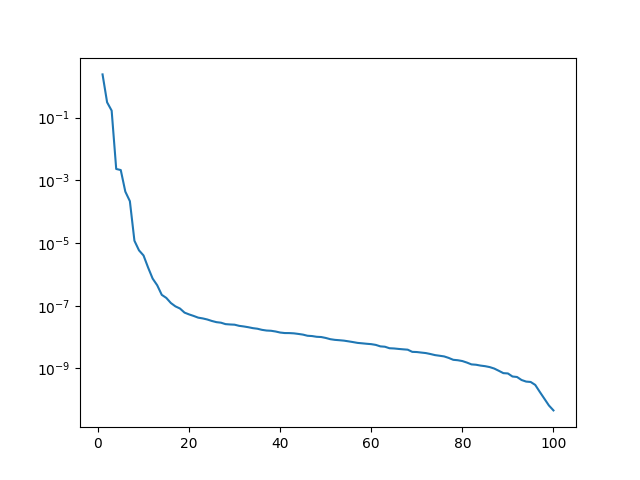}
  \subcaption{Singular value against singular value index.}
  \label{fig:1dpodsing}
\end{minipage}%
\hfill
\begin{minipage}{.45\textwidth}
  \centering
  \includegraphics[scale=0.5]{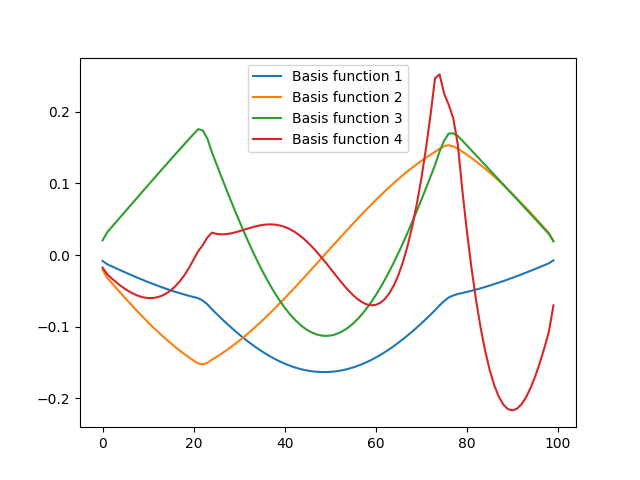}
  \subcaption{First four basis functions.}
  \label{fig:1dpodbasis}
\end{minipage}
\caption{Singular values and basis functions of snapshots from the 1D slab reactor.}
\label{1dpodbas}
\end{figure}

An SVD is applied to the 100 snapshots and from a possible 100 basis functions, 10 are retained, which corresponds to capturing 99.999\% of the information contained in the snapshots, see Figure~\ref{fig:1dpodsing}. As can be seen from this figure, the magnitude of the singular values decreases rapidly from the 1st to the 20th singular value (by about 7 orders of magnitude) before levelling off, and it can be expected from this that a POD model with a sufficient number of basis functions will be perform well. The first four basis functions can be seen in Figure~\ref{fig:1dpodbasis}. Notice that the fourth basis function has more detail or structure than the other basis functions, following the general trend that the higher-order basis functions tend to have more structure.  

\begin{figure}[!htb]
    \centering
    \includegraphics[scale=0.5]{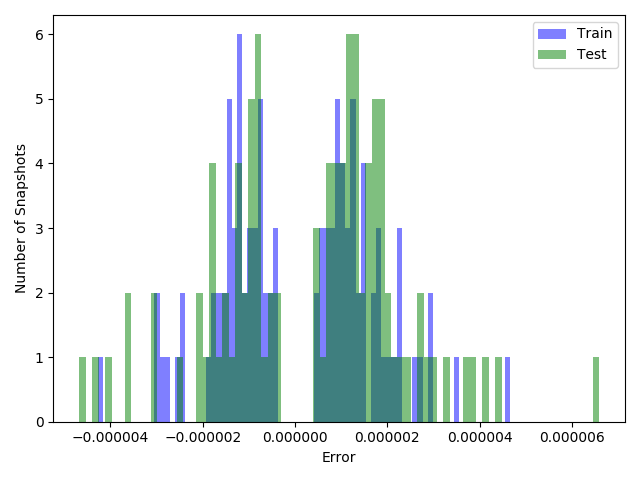}
    \caption{The error introduced to the representation of the high-fidelity model solutions by projecting onto the reduced space defined by the POD basis functions for the 1D test case. Seen data (snapshots) in blue, unseen data in green.}
 \label{fig:1dpod_comp}
\end{figure}

Figure~\ref{fig:1dpod_comp} shows the error, $e^{\text{POD}}_{\text{max}}(\bm{\phi}^{\text{HFM}})$,  resulting from projecting high-fidelity model solutions onto the reduced space for both the seen data (snapshots or training data) and the unseen data (test data), see Equation~\eqref{error_compression_POD}. Both data sets have a similar range of errors and the errors themselves are low:  $\mathcal{O}(\num{e-6})$.

The 200 sets of parameters that generated the high-fidelity model results are used again to obtain the matrices $\bm{A}$ and $\bm{B}$. The system of equations describing the POD-based ROM, given in~\eqref{pod_eig}, is solved for each problem using the power method described in Algorithms~\ref{alg:PM_outer} and~\ref{alg:PODROM_inner}, and Section~\ref{sec:PODGROM} using numpy's Gaussian elimination function to solve for the reduced variables. Half of the results have been seen by the model and half have not been seen.

\begin{figure}[!htb]
\centering
\begin{minipage}{.45\textwidth}
  \centering
  \includegraphics[scale=0.5]{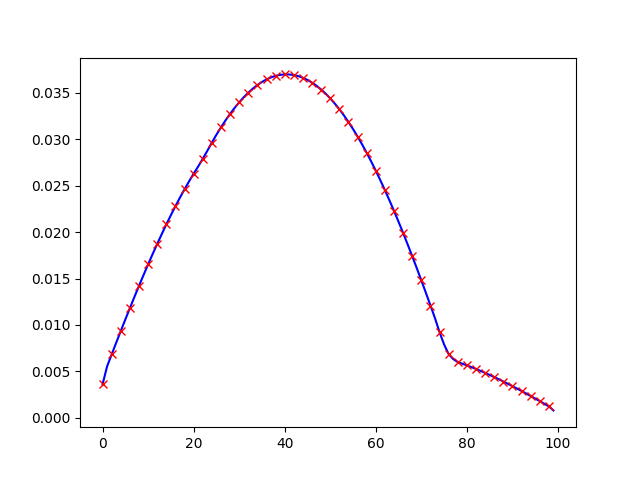}
  \subcaption{Scalar flux (\SI{}{neutrons.cm^{-2}.s^{-1}}) vs $x$ (\SI{}{\cm}), HFM solution in blue, POD-based ROM in red crosses.}
  \label{fig:dummyleft1}
\end{minipage}%
\hfill
\begin{minipage}{.45\textwidth}
  \centering
  \includegraphics[scale=0.5]{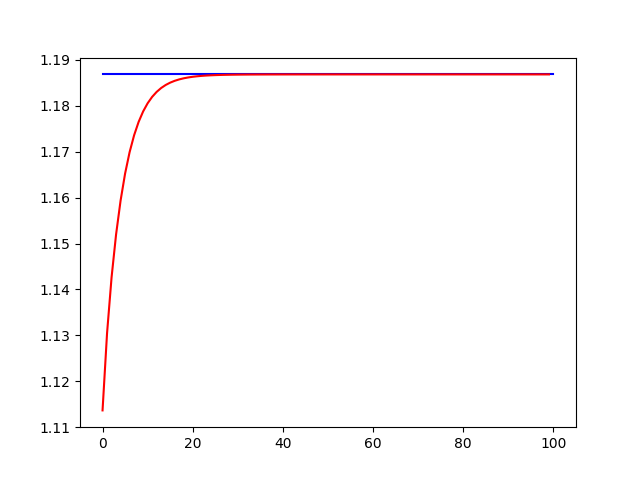}
  \subcaption{$k_{\text{eff}}$ vs iteration number, convergence of the POD-based ROM (red) to the converged HFM value (blue).}
  \label{fig:dummyright1}
\end{minipage}
\caption{Scalar flux and convergence of $k_{\text{eff}}$ for a seen problem with material parameters (cross-sections) determined by $r_1=0.957$ and $r_2=0.115$, comparing the POD-based ROM with the high-fidelity model.}
\label{pod30}
\end{figure}

\begin{figure}[!htb]
\centering
\begin{minipage}{.45\textwidth}
  \centering
  \includegraphics[scale=0.5]{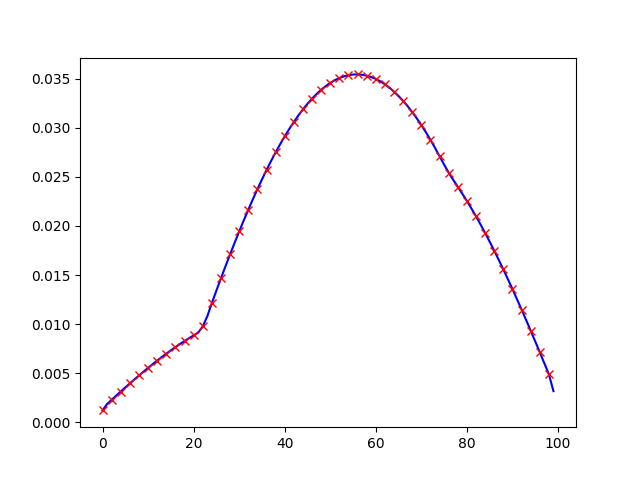}
  \subcaption{Scalar flux (\SI{}{neutrons.cm^{-2}.s^{-1}}) vs $x$ (\SI{}{\cm}), HFM solution in blue, POD-based ROM model in red crosses.}
  \label{dummy_seen1}
\end{minipage}%
\hfill
\begin{minipage}{.45\textwidth}
  \centering
  \includegraphics[scale=0.5]{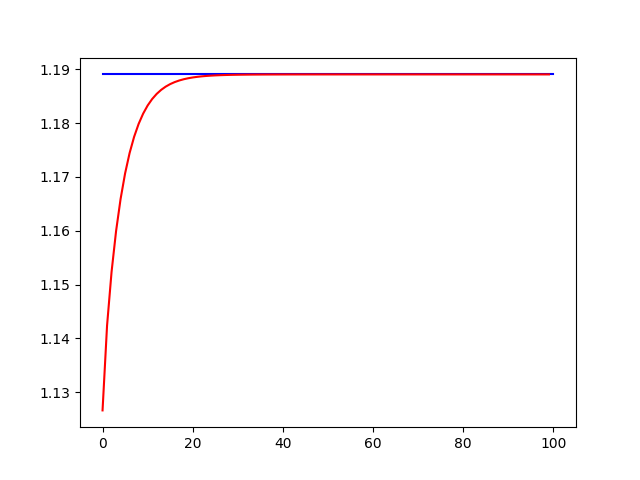}
  \subcaption{$k_{\text{eff}}$ vs iteration number, convergence of the POD-based ROM (in red) to the converged HFM value (in blue).}
  \label{dummy_pred1}
\end{minipage}
\caption{Scalar flux and convergence of $k_{\text{eff}}$ for an unseen problem with material parameters determined by $r_1=0.458$ and $r_2=0.932$ comparing the POD-based ROM with the high-fidelity model.}
\label{pod25}
\end{figure}

Figure~\ref{pod30} shows the flux profile and the convergence of $k_{\text{eff}}$ for the mixing coefficients $r_1=0.957$ and $r_2=0.115$. In this figure, results from the POD-based reduced-order model are compared with the high-fidelity model for this particular set of parameter values. This solution (of the high-fidelity model) was included in the snapshots, i.e.~is seen. Figure~\ref{pod25} shows the flux profile and the convergence of $k_{\text{eff}}$ for the material parameters $r_1=0.458$ and $r_2=0.932$, again, comparing the POD-based reduced-order model with the high-fidelity model. This solution (of the high-fidelity model) was not included in the snapshots, i.e.~is unseen. Both examples demonstrate excellent agreement of the POD-based scalar flux profiles with those of the high-fidelity model and show good convergence of $k_{\text{eff}}$ towards the converged $k_{\text{eff}}$ values of the high-fidelity model.

\begin{figure}[!htb]
\centering
\begin{minipage}{.45\textwidth}
  \centering
  \includegraphics[scale=0.5]{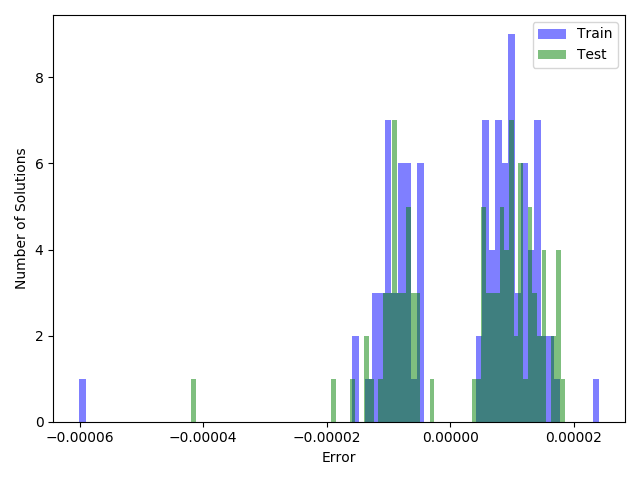}
  \subcaption{Histogram of the error in the flux profile.}
  \label{fig:1dpodint}
\end{minipage}%
\hfill
\begin{minipage}{.45\textwidth}
  \centering
  \includegraphics[scale=0.5]{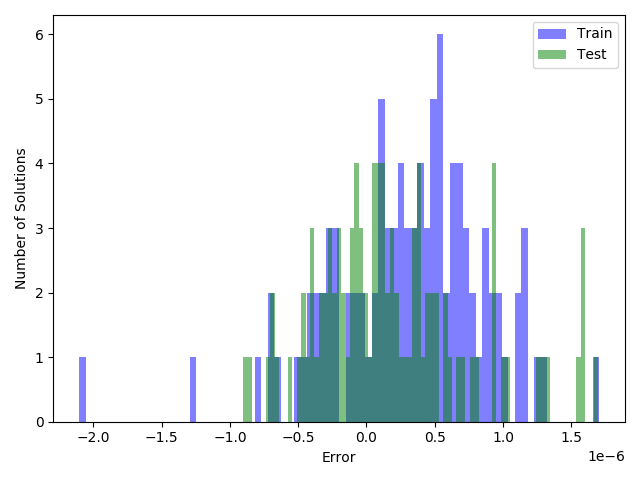}
  \subcaption{Histogram of the error in $k_{\text{eff}}$.}
  \label{fig:1dpodkeff}
\end{minipage}
\caption{Errors in the seen (blue) and unseen (green) results for the POD-based reduced-order model.}
\label{1dpod}
\end{figure}

To further analyze the results,  the solutions from the POD-based ROM were compared with solutions from the high-fidelity model (both seen and unseen data sets) using the error measures for the scalar flux and $k_{\text{eff}}$ as given by Equations~\eqref{error} and~\eqref{keff_error} respectively. Figure~\ref{fig:1dpodint} shows a histogram of the error in the scalar flux profile  and Figure~\ref{fig:1dpodkeff} shows a histogram of the error in $k_{\text{eff}}$. The plots  demonstrate both the accuracy of the POD-based model (for seen data) and its ability to prediction (for unseen parameters). The ranges of values of the errors for both the seen and unseen parameters are similar. Table~\ref{1dpod_err} displays the average maximum errors in the reconstruction, scalar flux and $k_{\text{eff}}$ as defined in Equations~\eqref{average_maximum_error}, \eqref{average_maximum_error_keff} and~\eqref{average_maximum_error_recon}. The projection error (for one solution) is given by Equation~\eqref{error_compression_POD}. The error in the scalar flux solution will contain this error and other numerical errors made when solving the reduced generalised eigenvalue problem. Therefore we expect that the projection error would be lower than the flux error, borne out by the results shown in Table~\ref{1dpod_err}. Also noteworthy is the similarity between the range of errors for the seen and unseen data. This may mean that the snapshots represent well the behaviour that was tested by the 100 unseen problems. For $k_{\text{eff}}$, the error for the unseen results is slightly lower that than of the seen results. This is surprising but could be explained by several outliers and the averaging of the maximum errors.

\begin{table}[h!]
\centering
 \begin{tabular}{l l c c} %{||l c c||} 
 \hline\hline
      & & Seen & Unseen  \\ [0.5ex] 
 \hline\hline
 Projection Error & $\bar{e}_{\text{max}}^{\text{POD}}(\bm{\phi}^{\text{HFM}})$ & \num{1.4344e-6} & \num{1.7004e-6}  \\[2mm] 
 %\hline
 Flux Error & $\bar{e}_{\text{max}}(\bm{\phi}^{\text{POD}})$ & \num{1.0270e-05} & \num{1.0401e-05}  \\[2mm] 
 %\hline
 $k_{\text{eff}}$ Error & $\bar{e}(k_{\text{eff}}^{\text{POD}})$ & \num{5.2956e-07} & \num{4.4348e-07}  \\[2mm] 
 \hline\hline
\end{tabular}
\caption{Average maximum errors for seen and unseen data using the POD-based reduced-order model.%to reduce 1D case from 100 variables to 10 variables.
}
\label{1dpod_err}
\end{table}

\subsubsection{Autoencoder-based reduced-order model}\label{1d_ae}

In this method the dimensionality reduction (or compression) is performed by an autoencoder instead of an SVD. Making this change requires an additional iteration loop when solving the reduced generalised eigenvalue problem as explained in Section~\ref{Proj_intr}. For the online prediction, an initial flux guess is compressed using the encoder. Following this, the $\bm{C}$ matrix is formed by using Equation~\eqref{eq:C}. The system in Equation~\eqref{eq:calculate_C} is then solved, the flux is decoded and passed to the outer iterations, and the whole process is repeated until $k_{\text{eff}}$ converges. In order to compare with results from the previous section, the autoencoder compresses to 10 latent variables. Autoencoders with different architectures were trained and the best performing architecture comprised 14 fully-connected layers with the following number of neurons in each layer: 
\begin{equation}
100 \rightarrow{} 100 \rightarrow{} 70 \rightarrow{} 50 \rightarrow{} 30 \rightarrow{} 20 \rightarrow{} 16 \rightarrow{} 10 \rightarrow{} 16 \rightarrow{} 20 \rightarrow{} 30 
\rightarrow{} 50 \rightarrow{} 70 
 \rightarrow{} 100 \rightarrow{} 100\ .
    \nonumber
\end{equation}
The activation function for every layer was the exponential linear unit (ELU)~\cite{clevert2015fast}; the optimiser used was `Nadam'~\cite{dozat2016incorporating} and the loss function was defined to be the mean squared reconstruction error. The snapshots were scaled between 0 and 1. The network was trained over 10,000 epochs using gradient descent (with a batch size of 100). In one epoch, all the training data is presented to the neural network.

\begin{figure}[!htb]
\centering
\begin{minipage}{.45\textwidth}
  \centering
  \includegraphics[scale=0.5]{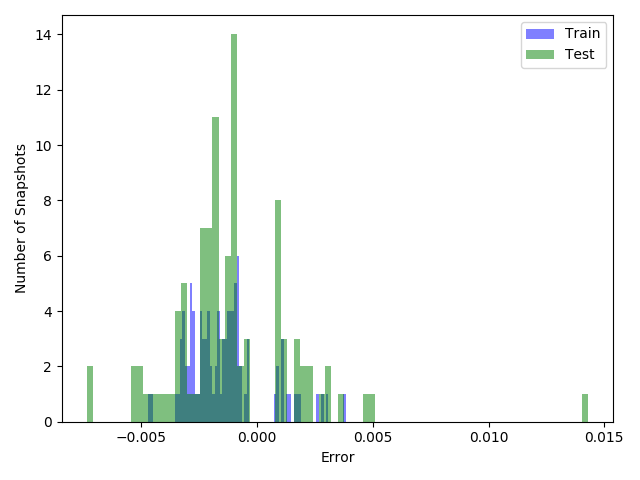}
  \subcaption{Compression error in the flux profile output by the autoencoder.}
  \label{fig:1daepass}
\end{minipage}%
\hfill
\begin{minipage}{.45\textwidth}
  \centering
  \includegraphics[scale=0.5]{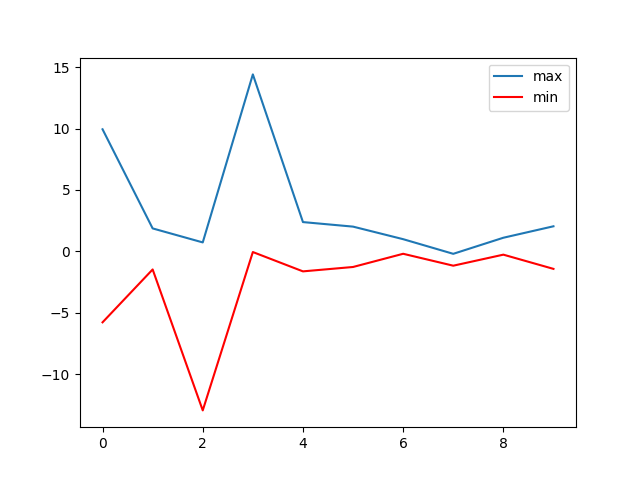}
  \subcaption{Minimum and maximum values of the 10 latent variables generated from the autoencoder over the training data.}
  \label{fig:1dlv}
\end{minipage}
%\caption{Verification and view at latent space of the standard autoencoder 1D-AE trained on 1D train dataset.}
\caption{Results from the autoencoder applied to the 1D test case.}
\label{1daever}
\end{figure}

A similar procedure is carried out as in the previous section and the same solutions of the high-fidelity model are used as seen data (the snapshots) and unseen data. When training neural networks, the seen data is usually referred to as training data and the unseen data is referred to as test data. The autoencoder was trained on the 100 snapshots and the other 100 unseen high-fidelity model results were used to test the prediction power of the AE-based reduced-order model. Figure~\ref{fig:1daepass} shows the error (as defined by Equation~\eqref{error_compression_AE}) between each HFM solution and the corresponding output of the autoencoder. The results are split into seen or training data and unseen or test data. The autoencoder performs well, achieving, at worst, an error of 1.5\% in the flux, with almost all of the absolute values of the error less than 0.5\% . Although the training and test data have a similar range of errors, it should be noted that these errors are much larger than those of the SVD (see figure~\ref{fig:1dpod_comp}). %, and the ROM based on this autoencoder will not be more accurate than this. 
Shown in Figure~\ref{fig:1dlv} are the minimum and maximum values taken by the 10~latent variables over the training data. Most of the latent variables here are contributing to the system although a couple of the variables have minimum and maximum values close to zero. 

\begin{figure}[!htb]
\centering
\begin{minipage}{.45\textwidth}
  \centering
  \includegraphics[scale=0.5]{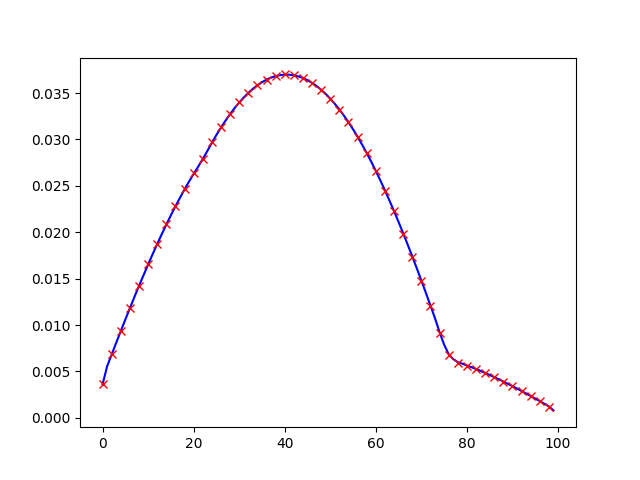}
  \subcaption{Scalar flux (\SI{}{neutrons.cm^{-2}.s^{-1}}) vs x(\SI{}{cm}). HFM solution in blue, AE-based ROM model in red crosses.}
  \label{fig:dummyleft2}
\end{minipage}%
\hfill
\begin{minipage}{.45\textwidth}
  \centering
  \includegraphics[scale=0.5]{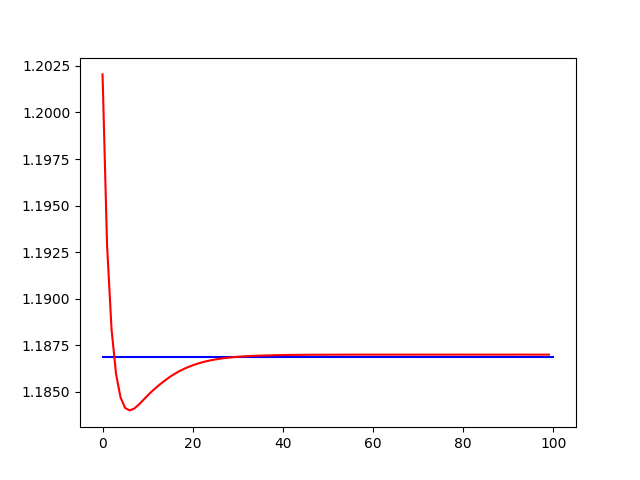}
  \subcaption{$k_{\text{eff}}$ vs iteration number, convergence of the AE-based ROM  (red) to the converged HFM value (blue).}
  \label{fig:dummyright2}
\end{minipage}
\caption{Flux and  convergence  of $k_{\text{eff}}$ for seen material parameters $r_1=0.957$ and $r_2=0.115$ for the AE-based reduced-order model and the high-fidelity model.}
\label{ae30}
\end{figure}

\begin{figure}[!htb]
\centering
\begin{minipage}{.45\textwidth}
  \centering
  \includegraphics[scale=0.5]{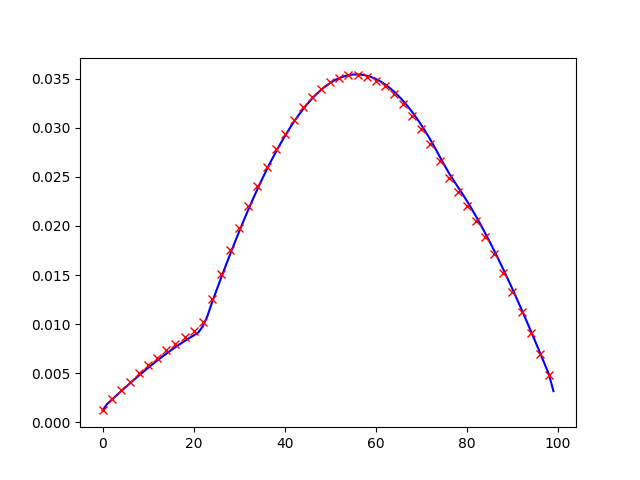}
  \subcaption{Scalar flux (\SI{}{neutrons.cm^{-2}.s^{-1}}) vs $x$ (\SI{}{\cm}). HFM solution in blue, AE-based ROM model in red crosses.}
  \label{dummy_seen2}
\end{minipage}%#
\hfill
\begin{minipage}{.45\textwidth}
  \centering
  \includegraphics[scale=0.5]{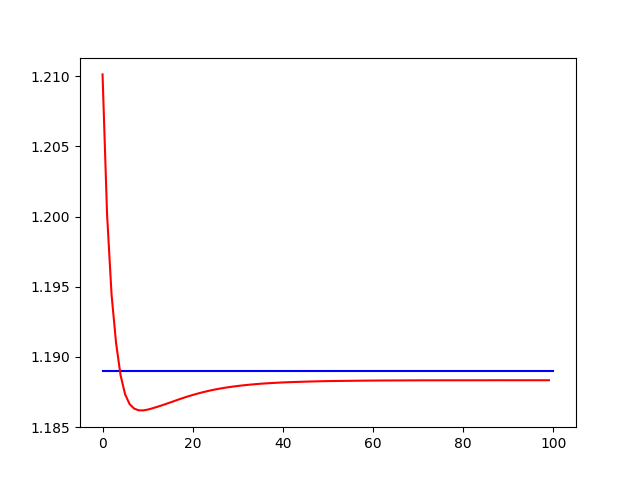}
  \subcaption{$k_{\text{eff}}$ vs iteration number, convergence of the AE-based ROM (red) to the converged HFM value (blue).}
  \label{dummy_pred2}
\end{minipage}
\caption{Scalar flux and convergence of $k_{\text{eff}}$ for an unseen problem with material parameters $r_1=0.458$ and $r_2=0.932$ comparing the AE-based reduced-order model and the high-fidelity model.}
\label{ae25}
\end{figure}

Figures~\ref{ae30} and~\ref{ae25} show the flux profile and the convergence of $k_{\text{eff}}$ for two sets of parameters, the first was used in the training (seen) and the second from the test data (unseen). Both seen and unseen flux profiles of the reduced-order model agree well with the high-fidelity model solution, although the convergence of the reduced-order model is very slightly worse for the unseen case, and the value of $k_{\text{eff}}$ is slightly less accurate.

\begin{figure}[!htb]
\centering
\begin{minipage}{.45\textwidth}
  \centering
  \includegraphics[scale=0.5]{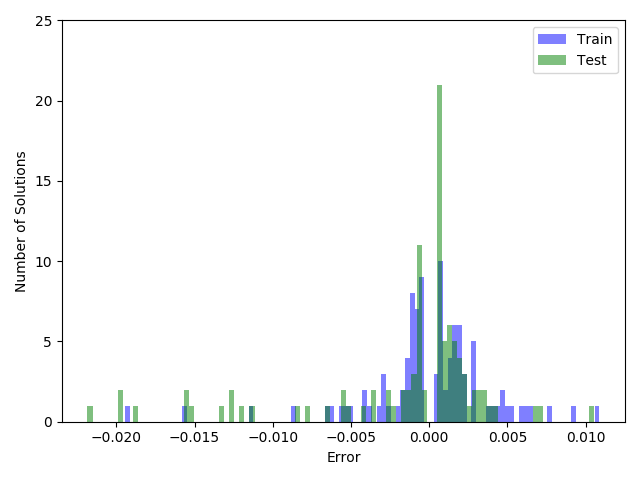}
  \subcaption{Histogram of the error in the flux profile.}
  \label{fig:1aeint}
\end{minipage}%
\hfill
\begin{minipage}{.45\textwidth}
  \centering
  \includegraphics[scale=0.5]{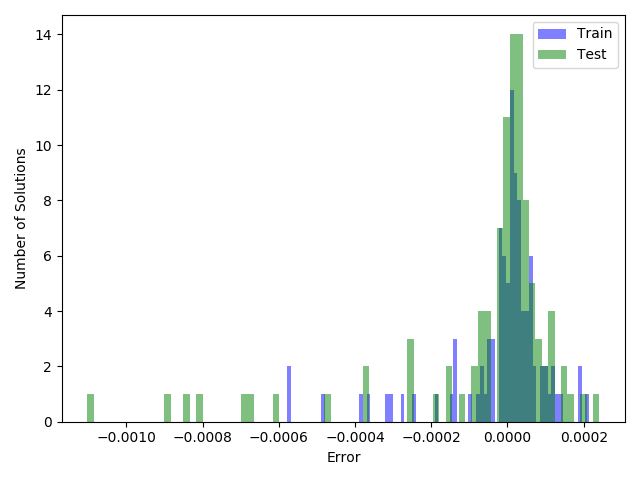}
  \subcaption{Histogram of the error in $k_{\text{eff}}$.}
  \label{fig:1daekeff}
\end{minipage}
\caption{Errors in the seen and unseen results for the autoencoder-based reduced-order model.}
%\caption{Errors for train and test data sets using intrusive method with the standard autoencoder 1D-AE.}
\label{AE_int}
\end{figure}

Figure~\ref{fig:1aeint} shows the error in flux and  Figure~\ref{fig:1daekeff} shows the error in $k_{\text{eff}}$ based on the error measures defined in Equations~\eqref{error} and~\eqref{keff_error}. 
Excluding outliers, most of the flux solutions have an absolute error of 1\% or less, and the absolute error in the $k_{\text{eff}}$ values is less than \num{0.4e-3}. These values are higher than for the POD-based reduced-order model, and there is also a larger spread of errors for the AE-based model than for the POD-based method. Table~\ref{1dae_err} shows the average maximum errors for the AE-based reduced-order model (see Equations~\eqref{average_maximum_error}, \eqref{average_maximum_error_keff} and~\eqref{average_maximum_error_recon}).  It can be seen that the scalar flux error is of the same magnitude as the compression error. Here, as for the POD-based ROM, the errors for the training data (seen) are very close to those of the test data (unseen) which, again, would seem to confirm that the training data captures most of the significant behaviour that we test with the unseen data set. Again, the compression error acts as a lower bound for the flux error. On the whole, these errors are one or two orders of magnitude larger than those seen for the POD-based ROM.

\begin{table}[h!]
\centering
 \begin{tabular}{l l c c} 
 \hline\hline
      &  & Seen/Training data & Unseen/Test data  \\ [0.5ex] 
 \hline\hline
 Compression Error & $\bar{e}^{\text{AE}}_{\text{max}}(\bm{\phi}^{\text{HFM}})$ &  \num{1.9143e-3} & \num{2.3166e-3}  \\[1mm]
 %\hline
  Flux Error & $\bar{e}_{\text{max}}(\bm{\phi}^{\text{AE}})$ & \num{2.7746e-3} & \num{3.6602e-3}  \\[1mm]
 %\hline
 $k_{\text{eff}}$ Error & $\bar{e}(k_{\text{eff}}^{\text{AE}})$ & \num{8.0772e-05} & \num{1.2149e-4}  \\[1mm]
 \hline\hline
\end{tabular}
\caption{Average maximum errors for seen and unseen data using the autoencoder-based reduced-order model.}
%\caption{Average maximum errors for train and test data sets using intrusive method with the standard autoencoder 1D-AE trained on 1D train dataset to reduce 1D case from 100 variables to 10 variables.}
\label{1dae_err}
\end{table}

For this test case, with two control-rod regions, one could postulate that having a low-dimensional space of dimension~2 should be enough to describe the behaviour of the system. The next section compares a projection-based ROM and an autoencder-based ROM with two basis functions and latent variables respectively.

\subsubsection{Reduction to two variables}\label{1d_2var} 
%The results of the variational autoencoder in Section~\ref{1d_vae} have indicated that two latent variables can be used to capture almost all the information contained in the training data. 
In order to see how well the models perform reducing to a low-dimensional space of dimension~2, a POD-based ROM and an AE-based ROM were constructed using just two POD basis functions for POD, and two latent variables for the autoencoder. 

Table~\ref{2var_err} shows the errors for compression, flux and  $k_{\text{eff}}$ when the system is compressed from 100 variables to two variables. For this extreme case, it can be seen that the autoencoder performs better than POD, and for $k_{\text{eff}}$, has errors that are an order of magnitude lower than for POD. This suggests that the autoencoder can capture more of the features than POD with a given number of basis functions/latent variables.
\begin{table}[h!]
\centering
 \begin{tabular}{l l c c l c c  } 
 \hline\hline
      &  \multicolumn{3}{c}{POD} & \multicolumn{3}{c}{AE}\\
 \hline\hline
     & & Seen & Unseen & & Train & Test  \\ 
 %\hline
 Compression Error & $\bar{e}^{\text{POD}}_{\text{max}}(\bm{\phi}^{\text{HFM}})$ & \num{6.3880e-2} & \num{7.1400e-2} & $\qquad\bar{e}^{\text{AE}}_{\text{max}}(\bm{\phi}^{\text{HFM}})$ & \num{4.0990e-3} & \num{4.2959e-3} \\[1mm] 
 %\hline
 Flux Error & $\bar{e}_{\text{max}}(\bm{\phi}^{\text{POD}})$ & \num{6.9400e-2} & \num{7.7714e-2} & $\qquad\bar{e}_{\text{max}}(\bm{\phi}^{\text{AE}})$ & \num{1.0946e-2} & \num{1.1368e-2}\\[1mm] 
 %\hline
 $k_{\text{eff}}$ Error & $\bar{e}(k_{\text{eff}}^{\text{POD}})$ & \num{9.8282e-3} & \num{1.0248e-2} & $\qquad\bar{e}(k_{\text{eff}}^{\text{AE}})$ & \num{8.6002e-4} & \num{8.2654e-4}\\[1mm]
 \hline\hline

\end{tabular}
\caption{Average maximum errors for seen results / training data and unseen results / test data using POD and a standard autoencoder to reduce from 100 variables to two variables.}
\label{2var_err}
\end{table}

\subsection{2D reactor core eigenvalue problem}\label{2d_case} 
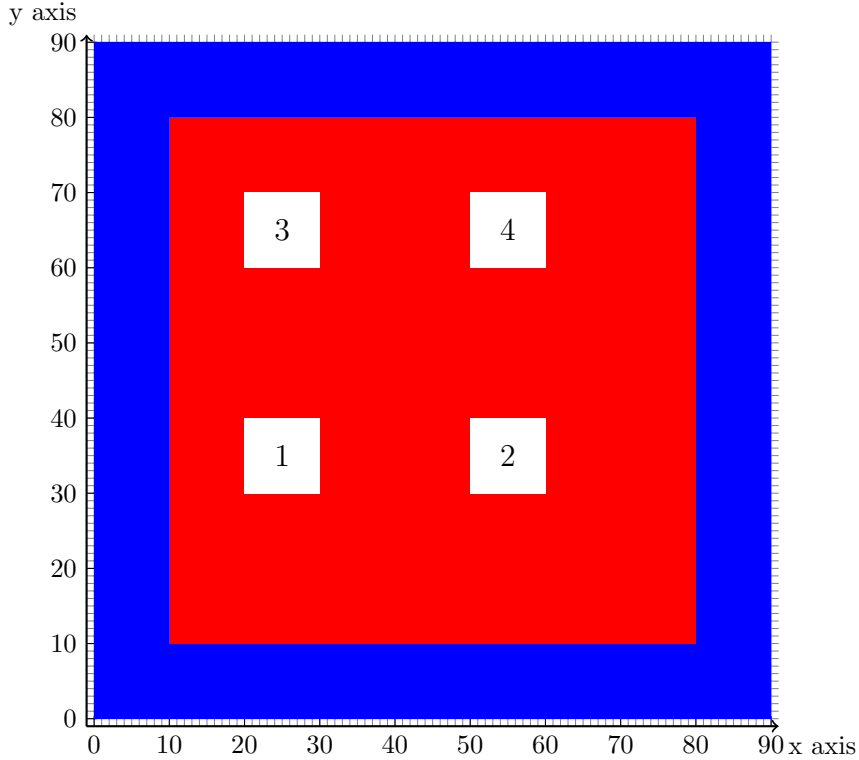
\begin{figure}[!htb]
\def\layersep{2.5cm}
\centering
\begin{tikzpicture}
\draw[step=0.1cm,gray,very thin] (3.9,3.9) grid (13.1,13.1);
\draw[thick,->] (3.9,3.9) -- (13.1,3.9) node[anchor=north west] {x axis};
\draw[thick,->] (3.9,3.9) -- (3.9,13.1) node[anchor=south east] {y axis};
 \foreach \x/\xtext in {4/0,5/10,6/20,7/30,8/40,9/50,10/60,11/70,12/80,13/90}
       \draw[yshift=4cm] (\x cm,1pt) -- (\x cm,-3pt) node[anchor=north] {$\xtext$};
   \foreach \y/\ytext in {4/0,5/10,6/20,7/30,8/40,9/50,10/60,11/70,12/80,13/90}
        \draw [xshift=4cm](1pt,\y cm) -- (-3pt,\y cm) node[anchor=east] {$\ytext$};
%\fill[pattern=north west lines, pattern color=black] (4,4) rectangle (13,13);
\fill[blue] (4,4) rectangle (13,13);
\fill[red] (5,5) rectangle (12,12);
\fill[white] (6,7) rectangle (7,8);
\fill[white] (6,10) rectangle (7,11);
\fill[white] (9,7) rectangle (10,8);
\fill[white] (9,10) rectangle (10,11);
\node at (6.5,7.5) {\large{1}};
\node at (6.5,10.5) {\large{3}};
\node at (9.5,7.5) {\large{2}};
\node at (9.5,10.5) {\large{4}};
\end{tikzpicture}
\caption{Geometry of a 2D reactor core with graphite (blue), fuel (red)  and four control-rod regions labelled 1,2,3,4 (white).}
\label{fig:2dcore}
\end{figure}

The methods are now applied to a simple 2D reactor. The domain is square-shaped with sides of length \SI{90}{cm} and there are four different materials in this mock reactor, see Figure~\ref{fig:2dcore}.  The domain is discretised with 90 cells in both directions, each cell is of length \SI{1}{cm}, making 8100~cells. In addition to this, 364 cells were used to enforce the boundary conditions. Once again, this test case was used by Buchan~et~al.~\cite{buchan2013pod} who based the cross-sections on IAEA benchmarks. The same macroscopic cross-sections are used here and their values are given in Table~\ref{2dprop}. 

\begin{table}[h!]
\centering
 \begin{tabular}{l l l l} 
 \hline\hline
         & $\Sigma_a$ & $\Sigma_s$ & $\Sigma_f$ \\ [0.5ex] 
 \hline\hline
 Fuel & 0.075 & 0.53 & 0.79 \\ 
 %\hline
 Water & 0.01 & 0.89 & 0 \\ 
 %\hline
 Control Rod & 0.38 & 0.2 & 0 \\ 
 %\hline
 Graphite & 0.15 & 0.5 & 0 \\[1ex] 
 \hline\hline
\end{tabular}
\caption{Macroscopic cross-sections for the materials in the 2D reactor core (units \si{\per\cm}).}
\label{2dprop}
\end{table}

In Figure~\ref{fig:2dcore} the blue region represents graphite, the red area contains the fuel and the four control-rod regions are white and labelled 1, 2, 3 and 4. Each control-rod region is \SI{10}{cm} square with centres at (25,35), (55,35), (25,65) and  (55,65) respectively. These are assumed to have a uniform distribution of the averaged properties of the mixture within the system based on a mixing coefficient for each region chosen. The method for determining the mixing coefficient is described in Section~\ref{sec:homogenisation}, however, instead of mixing fuel and control rod materials in the control-rod regions as done for the 1D test case, here, water and control rod are mixed, with more weight given to the water in order  to prevent domination of the strong absorption of the control rods.

For this system 800 high-fidelity model solutions are generated, each with four distinct values of~$r$, one for each control-rod region. The method for generating the solutions is similar to Section~\ref{1D_slab}. Of the 800 solutions, 400 are used during the offline phase as the snapshots (or seen or training data). The other 400 results make up the unseen data or test data.  These high-fidelity model solutions are used to construct a POD-based ROM and two autoencoder-based ROMs: the first uses a fully-connected autoencoder, the second uses the hybrid SVD-autoencoder. The number of degrees of freedom of the high-fidelity model solutions is 8100, and, for all ROMs, the number of POD coefficients or latent variables is chosen as 4, as this is perhaps the minimum number required to represent the system. 

\subsubsection{POD-based ROM}\label{2d_POD} 
For the POD-based ROM, an SVD is used to compress the snapshots. From a possible 400 basis functions, 4 are retained which corresponds to capturing 99.878\% of the information contained in the snapshots. In the online stage the parameter values are chosen and matrices $\bm{A}$ and $\bm{B}$ of the high-fidelity model are formed. The first 4 basis functions are used to create $\bm{R}$ and Equation~\eqref{pod_eig} is formed. This is solved with the power method as described in Algorithm~\ref{alg:PODROM_inner} and Section~\ref{ROM-POD}. 

\begin{figure}[!htb]
\centering
\begin{minipage}{.45\textwidth}
  \centering
  \includegraphics[scale=0.5]{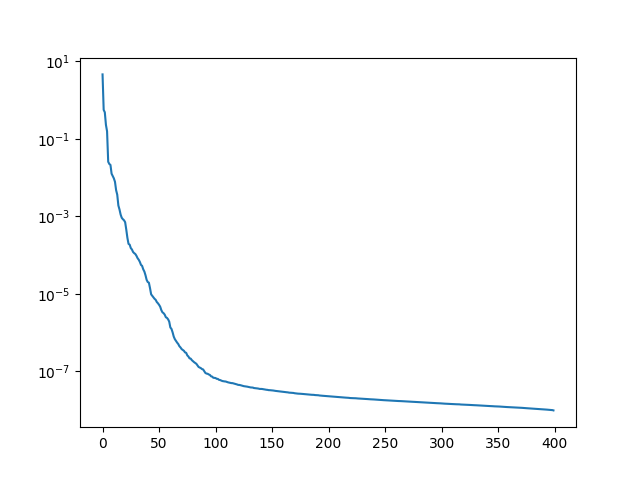}
  \subcaption{Singular value against singular value index.}
  \label{fig:2dpodsing}
\end{minipage}%
\hfill
\begin{minipage}{.45\textwidth}
  \centering
  \includegraphics[scale=0.5]{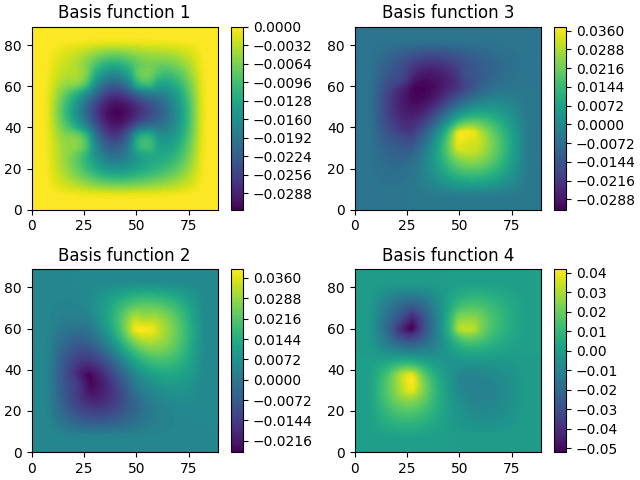}
  \subcaption{First four basis functions.}
  \label{fig:2dpodbasis}
\end{minipage}
\caption{Singular values and four basis functions for the 2D test case.}
\label{2dpodbas}
\end{figure}

\begin{figure}[!htb]
    \centering
    \includegraphics[scale=0.5]{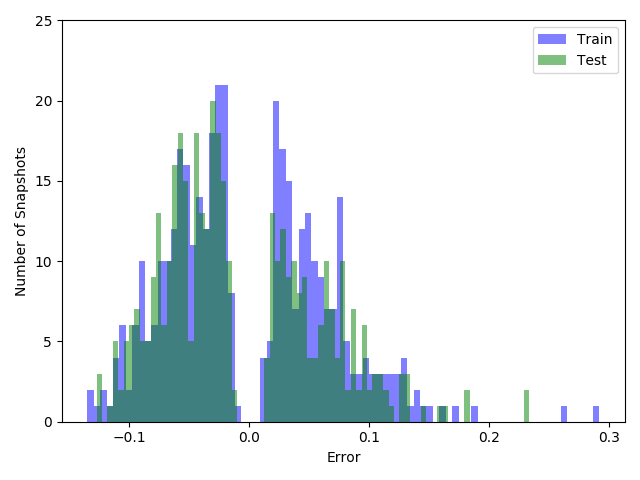}
    \caption{The error introduced to the representation of the high-fidelity model solutions by projecting onto the reduced space defined by the POD basis functions for the 2D problem. Seen data (snapshots) in blue, unseen data in green.}
 \label{fig:2dpod_comp}
\end{figure}

Figure~\ref{fig:2dpodsing} shows the singular values, of which the first 100 decrease rapidly, the remaining values reach a plateau. The first four corresponding basis functions can be seen in Figure~\ref{fig:2dpodbasis}. Figure~\ref{fig:2dpod_comp} shows the reconstruction error in truncating the POD basis for the snapshots (seen data) and the unseen data, see Equation~\eqref{error_compression_POD}. The range of errors for both sets of data is very similar.  

\begin{figure}[!htb]
\centering
\begin{minipage}{.45\textwidth}
  \centering
  \includegraphics[scale=0.5]{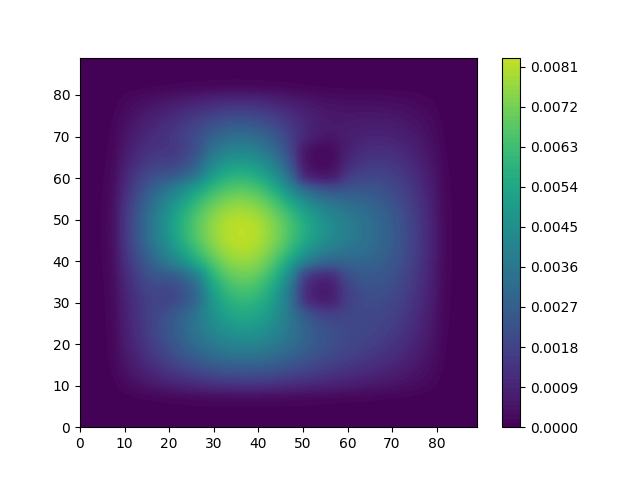}
  \subcaption{High-fidelity model.}
  \label{fig:dummyleft4}
\end{minipage}%
\hfill
\begin{minipage}{.45\textwidth}
  \centering
  \includegraphics[scale=0.5]{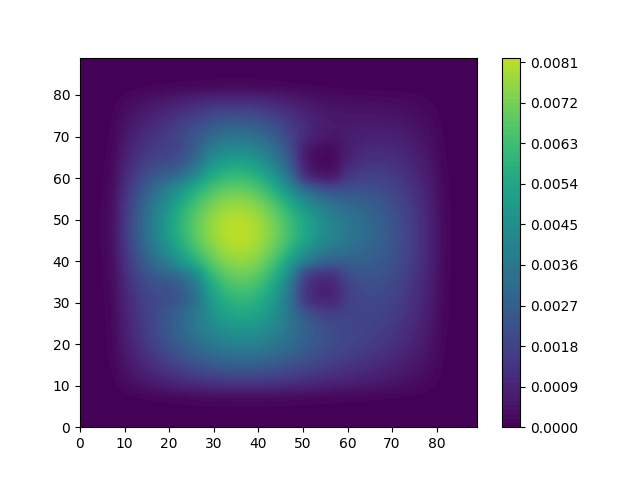}
  \subcaption{POD-based ROM.}
  \label{fig:dummyright4}
\end{minipage}
\caption{Scalar flux (\si{neutrons.cm^{-2}.s^{-1}}) across the reactor for a seen problem with material properties $r_1=0.9782$, $r_2=0.9891$, $r_3=0.7006$ and $r_4=0.8316$ using the POD-based ROM.}
\label{2dpodseen}
\end{figure}

\begin{figure}[!htb]
\centering
\begin{minipage}{.45\textwidth}
  \centering
  \includegraphics[scale=0.5]{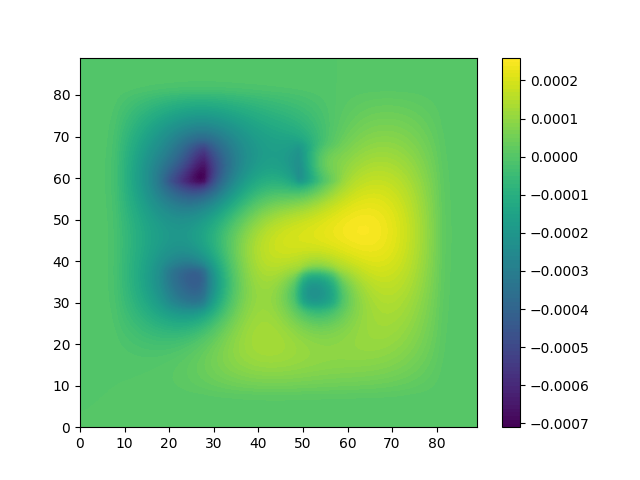}
  \subcaption{Pointwise error in scalar flux.}
  \label{fig:2dpod_diff_seen}
\end{minipage}%
\hfill
\begin{minipage}{.45\textwidth}
  \centering
  \includegraphics[scale=0.5]{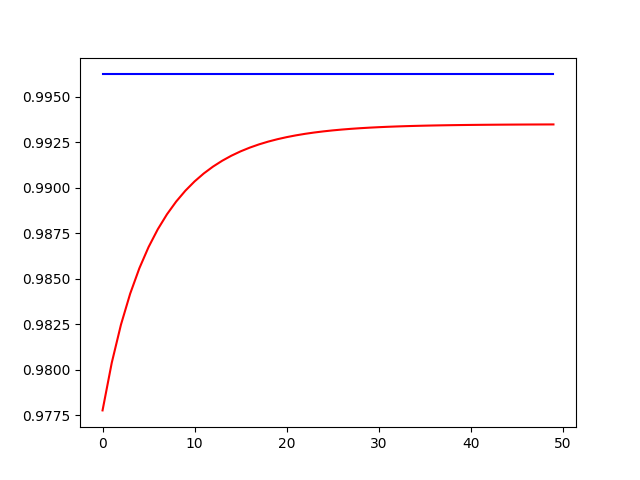}
  \subcaption{$k_{\text{eff}}$ vs iteration number, convergence of POD-based ROM (red) to the converged HFM value (blue).}
  \label{fig:2dpod_keff_seen}
\end{minipage}
\caption{Pointwise error in scalar flux and convergence of $k_{\text{eff}}$ for the POD-based ROM with material properties $r_1=0.9782$, $r_2=0.9891$, $r_3=0.7006$ and $r_4=0.8316$.}
\label{pod_keffconv_seen}
\end{figure}

The two flux profiles shown in Figure~\ref{2dpodseen}, for a seen set of parameter values, are in close agreement. The high-fidelity model results are on the left, results from the POD-based ROM are on the right. The pointwise error in the flux is shown in Figure~\ref{fig:2dpod_diff_seen}, and most error can be seen near the control-rod regions. The value of $k_{\text{eff}}$ converges to a value within 0.0035 of the value given by the high-fidelity model, see Figure~\ref{fig:2dpod_keff_seen}.


\begin{figure}[!htb]
\centering
\begin{minipage}{.45\textwidth}
  \centering
  \includegraphics[scale=0.5]{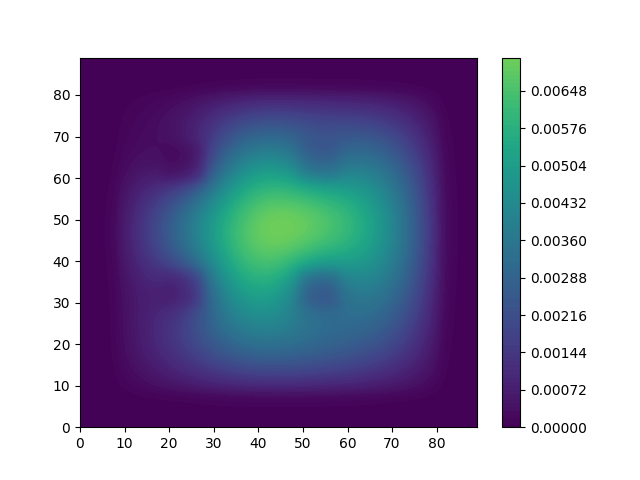}
  \subcaption{High-fidelity model.}
  \label{fig:dummyleft5}
\end{minipage}%
\hfill
\begin{minipage}{.45\textwidth}
  \centering
  \includegraphics[scale=0.5]{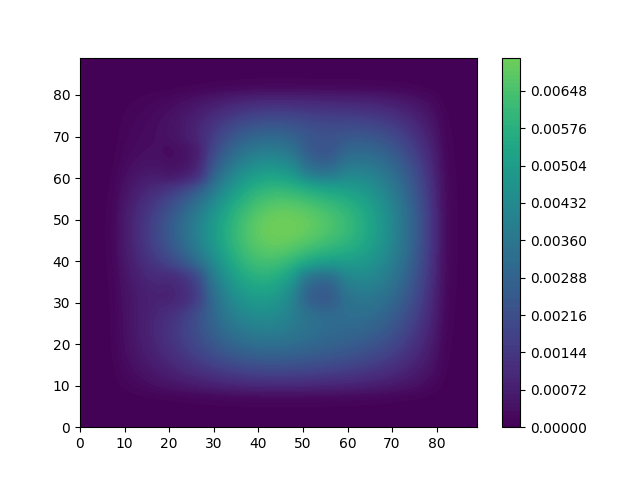}
  \subcaption{POD-based ROM.}
  \label{fig:dummyright5}
\end{minipage}
\caption{Scalar flux (\si{neutrons.cm^{-2}.s^{-1}}) across the reactor for an unseen problem with material properties $r_1=0.9452$, $r_2=0.8647$, $r_3=0.9996$ and $r_4=0.9776$ using POD-based ROM.}
\label{2dpodunseen}
\end{figure}

\begin{figure}[!htb]
\centering
\begin{minipage}{.45\textwidth}
  \centering
  \includegraphics[scale=0.5]{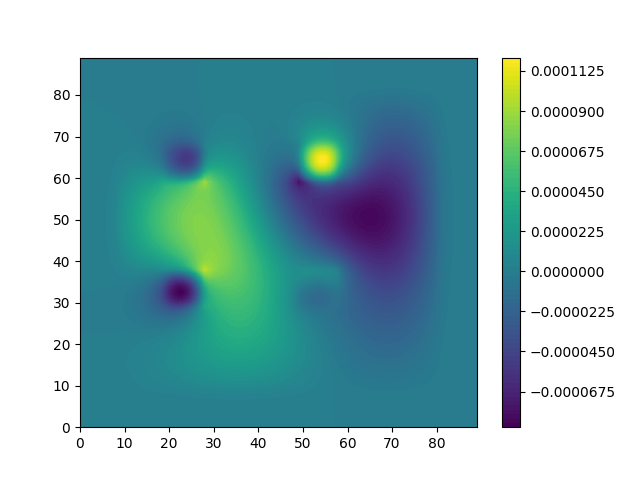}
  \subcaption{Pointwise error in scalar flux across the reactor.}
  \label{fig:2dpod_diff_unseen}
\end{minipage}%
\hfill
\begin{minipage}{.45\textwidth}
  \centering
  \includegraphics[scale=0.5]{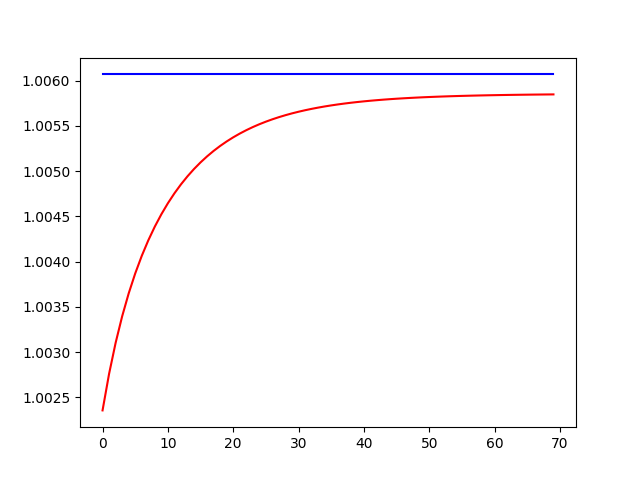}
  \subcaption{$k_{\text{eff}}$ vs iteration number, convergence of POD-based ROM (red) to the converged HFM value (blue).}
  \label{fig:2dpod_keff_unseen}
\end{minipage}
\caption{Pointwise error in scalar flux and convergence of $k_{\text{eff}}$ for an unseen problem with material properties $r_1=0.9452$, $r_2=0.8647$, $r_3=0.9996$ and $r_4=0.9776$ using POD-based ROM.}
\label{pod_conv_keff_unseen}
\end{figure}

Figure~\ref{2dpodunseen} shows flux profiles from the high-fidelity model and the POD-based ROM for an unseen set of parameters. The flux profiles are again close, with the pointwise error shown in Figure~\ref{fig:2dpod_diff_unseen}, where most of the error is concentrated around the control-rod regions. In  Figure~\ref{fig:2dpod_keff_unseen}, $k_{\text{eff}}$ is seen to converge closely to the value obtained from the high-fidelity model.


\begin{figure}[!htb]
\centering
\begin{minipage}{.45\textwidth}
  \centering
  \includegraphics[scale=0.5]{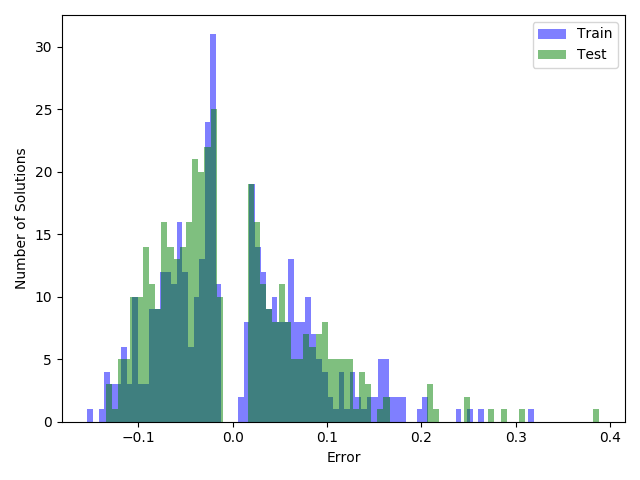}
  \subcaption{Histogram of the error in the flux profile.}
  \label{fig:2dpodint}
\end{minipage}%
\hfill
\begin{minipage}{.45\textwidth}
  \centering
  \includegraphics[scale=0.5]{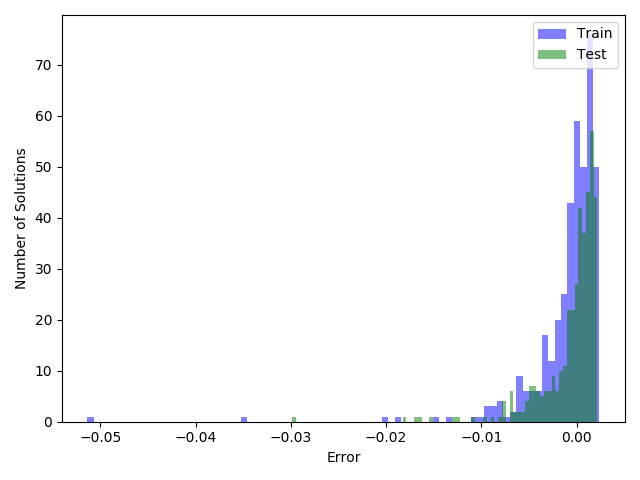}
  \subcaption{Histogram of the error in $k_{\text{eff}}$.}
  \label{fig:2dpodkeff}
\end{minipage}
\caption{Error in the seen (blue) and the unseen (green) results for the POD-based ROM applied to the 2D test case.}
\label{2dpod}
\end{figure}

Figure~\ref{2dpod} shows histogram plots of the errors from comparing the POD-based ROM with both the 400 snapshots (training data) and the 400 unseen HFM solutions (the test data). Figure~\ref{fig:2dpodint} shows the error in flux profile calculated using Equation~\eqref{error}, and Figure~\ref{fig:2dpodkeff} shows the error in $k_{\text{eff}}$ calculated using Equation~\eqref{keff_error}. Setting aside some outliers, the error in flux profile has a range similar to the range of errors seen in the truncation of the POD basis, implying that the main contribution to the error for the POD-based ROM is due to the truncation of the POD basis.

\subsubsection{Autoencoder-based reduced order model}\label{2d_ae} 
In the autoencoder-based ROM, the SVD is replaced by an autoencoder, which, for this case, compresses down from 8100~variables to 4~latent variables. The autoencoder comprises 14~fully connected layers with the following number of neurons in each layer: 
\begin{equation}
8100 \rightarrow{} 100 \rightarrow{} 70 \rightarrow{} 50 \rightarrow{} 30 \rightarrow{} 16 \rightarrow{} 8 \rightarrow{} 4 \rightarrow{} 8 \rightarrow{} 16 \rightarrow{} 30 
\rightarrow{} 50 \rightarrow{} 70 
 \rightarrow{} 100 \rightarrow{} 8100
    \nonumber
\end{equation}
Note here, the large reduction in neurons between the input and the first hidden layer. Although needed to reduce the number of layers which would otherwise be too high, this is expected to cause some difficulties with the accuracy of the network's predictions. The activation function for every layer is the exponential linear unit; the optimiser used is `Nadam' and the loss function is defined to be the mean squared reconstruction error. The snapshots are scaled between 0 and 1, and the network was trained using mini-batch gradient descent with a batch size of 50. The number of epochs was 30,000.

\begin{figure}[!htb]
\centering
\begin{minipage}{.45\textwidth}
  \centering
  \includegraphics[scale=0.5]{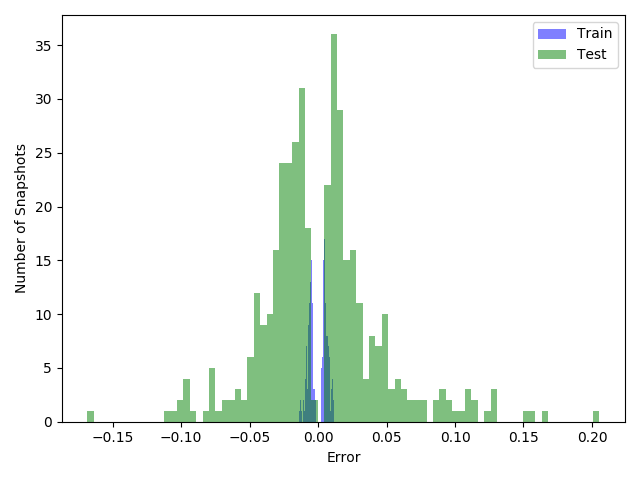}
  \subcaption{Reconstruction error in the flux profile for seen (blue) and unseen data (green).}
  \label{fig:2daepass}
\end{minipage}%
\hfill
\begin{minipage}{.45\textwidth}
  \centering
  \includegraphics[scale=0.5]{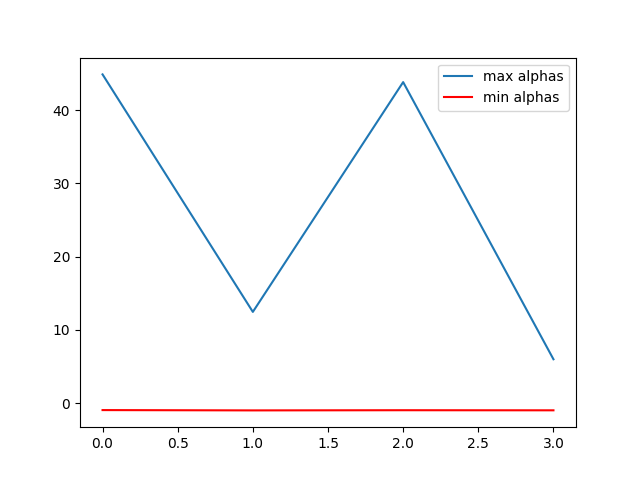}
  \subcaption{Minimum and maximum of the 4 latent variables generated from the autoencoder.}
  \label{fig:2daevat}
\end{minipage}
\caption{The representation error introduced to the representation of the high-fidelity model solutions when using the autoencoder, and the range of the latent variables.}
\label{ae4pass}
\end{figure}

Figure~\ref{fig:2daepass} shows the reconstruction or compression error (i.e.~the difference between each high-fidelity model solution and the corresponding output of the autoencoder after training) calculated using Equation~\eqref{error_compression_AE}. The range of errors in the test data is much larger (approximately ten times larger) than the range of errors for the training data, revealing that over-fitting has occurred. Figure~\ref{fig:2daevat} shows the minimum and maximum values attained by the latent variables over the training data, which shows the relative sizes and ranges of the latent variables.

\begin{figure}[!htb]
\centering
\begin{minipage}{.45\textwidth}
  \centering
  \includegraphics[scale=0.5]{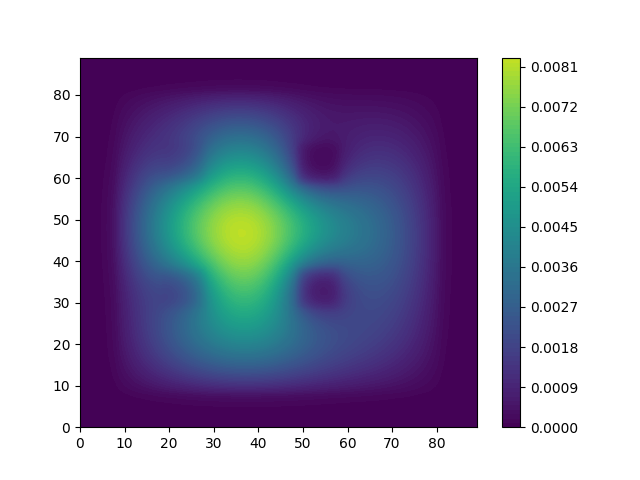}
  \subcaption{High-fidelity model.}
  \label{fig:dummyleft6}
\end{minipage}%
\hfill
\begin{minipage}{.45\textwidth}
  \centering
  \includegraphics[scale=0.5]{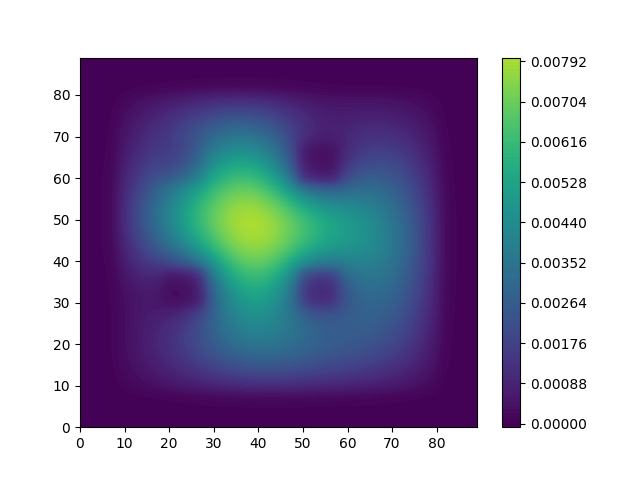}
  \subcaption{AE-based ROM.}
  \label{fig:dummyright6}
\end{minipage}
\caption{Scalar flux (\si{neutrons.cm^{-2}.s^{-1}}) across the reactor for a seen problem with material properties $r_1=0.9782$, $r_2=0.9891$, $r_3=0.7006$ and $r_4=0.8316$.}
\label{2daeseen}
\end{figure}

\begin{figure}[!htb]
\centering
\begin{minipage}{.45\textwidth}
  \centering
  \includegraphics[scale=0.5]{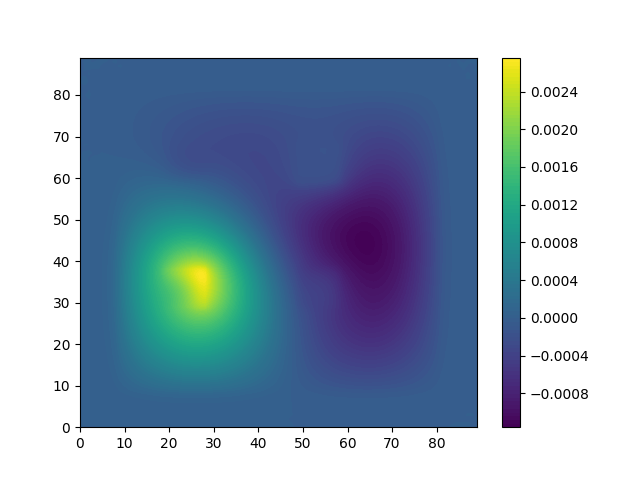}
  \subcaption{Pointwise error in scalar flux.}
  \label{fig:2dae_diff_seen}
\end{minipage}%
\hfill
\begin{minipage}{.45\textwidth}
  \centering
  \includegraphics[scale=0.5]{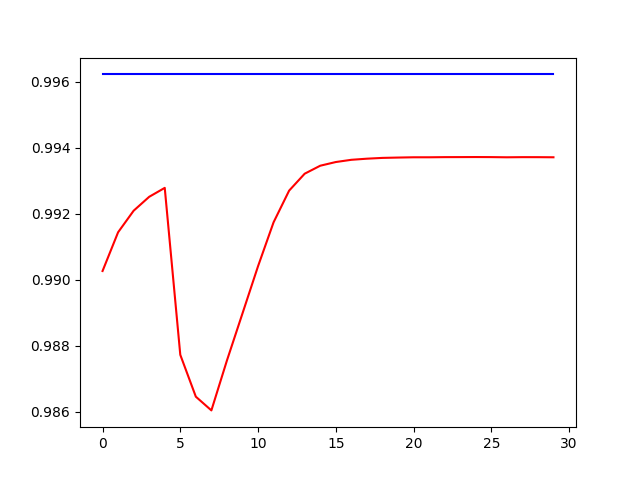}
  \subcaption{$k_{\text{eff}}$ vs iteration number, convergence of AE-based ROM (red) to the converged HFM value (blue).}
  \label{fig:2dae_keff_seen}
\end{minipage}
\caption{Pointwise error in scalar flux and convergence of $k_{\text{eff}}$ for the AE-based ROM with material properties $r_1=0.9782$, $r_2=0.9891$, $r_3=0.7006$ and $r_4=0.8316$.}
\label{ae_keffconv_seen}
\end{figure}

Results for a seen set of parameter values can be seen in Figures~\ref{2daeseen} and~\ref{ae_keffconv_seen}. There is a good agreement of the scalar flux profiles of the high-fidelity model and the AE-based ROM, with most error occurring in the vicinity of the control-rod regions (especially the lower two regions). This error is much more diffuse than for the POD-based ROM. The value of $k_{\text{eff}}$ seen in Figure~\ref{fig:2dae_keff_seen}, is in close agreement with that of the high-fidelity model.

\begin{figure}[!htb]
\centering
\begin{minipage}{.45\textwidth}
  \centering
  \includegraphics[scale=0.5]{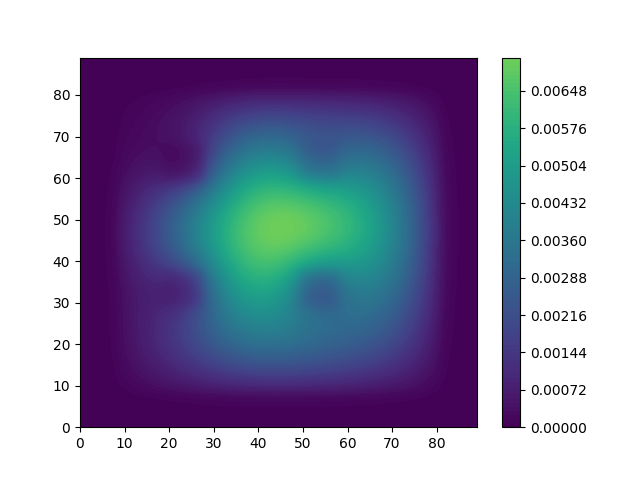}
  \subcaption{High-fidelity model.}
  \label{fig:dummyleft7}
\end{minipage}%
\hfill
\begin{minipage}{.45\textwidth}
  \centering
  \includegraphics[scale=0.5]{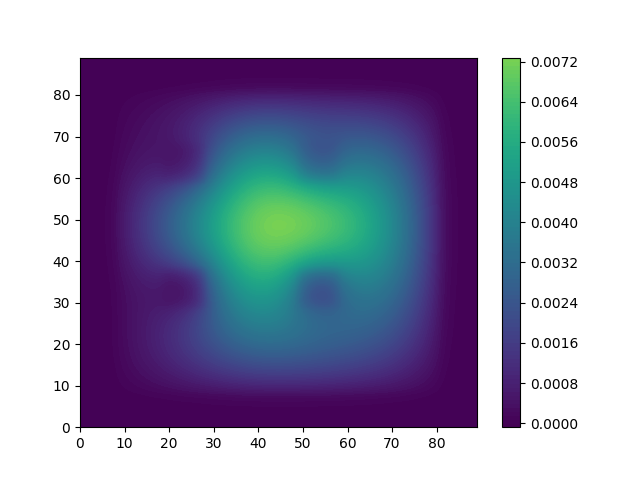}
  \subcaption{AE-based ROM.}
  \label{fig:dummyright7}
\end{minipage}
\caption{Scalar flux (\si{neutrons.cm^{-2}.s^{-1}}) across the reactor for an unseen problem with the AE-based ROM for material properties $r_1=0.9452$, $r_2=0.8647$, $r_3=0.9996$ and $r_4=0.9776$.}
\label{2daeunseen}
\end{figure}

\begin{figure}[!htb]
\centering
\begin{minipage}{.45\textwidth}
  \centering
  \includegraphics[scale=0.5]{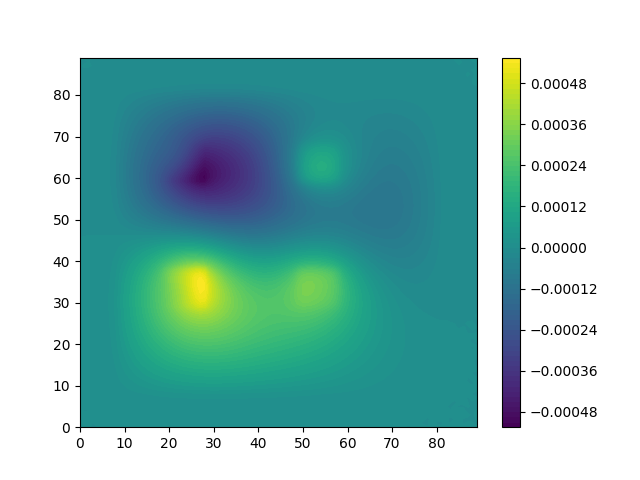}
  \subcaption{Pointwise error in scalar flux across the reactor.}
  \label{fig:2dae_diff_unseen}
\end{minipage}%
\hfill
\begin{minipage}{.45\textwidth}
  \centering
  \includegraphics[scale=0.5]{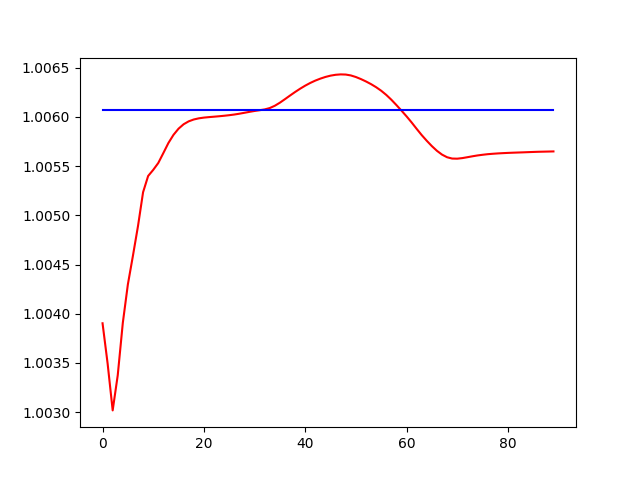}
  \subcaption{$k_{\text{eff}}$ vs iteration number, convergence of POD-based ROM (red) to the converged HFM value (blue).}
  \label{fig:2dae_keff_unseen}
\end{minipage}
\caption{Pointwise error in scalar flux and convergence of $k_{\text{eff}}$ for an unseen problem using the AE-based ROM with material properties $r_1=0.9452$, $r_2=0.8647$, $r_3=0.9996$ and $r_4=0.9776$.}
\label{2dae_conv_keff_unseen}
\end{figure}

The results for an unseen set of parameters can be seen in Figures~\ref{2daeunseen} and~\ref{2dae_conv_keff_unseen}. The pointwise error between the high-fidelity model and the AE-based ROM is shown in Figure~\ref{fig:2dae_diff_unseen} is again more diffuse than the error for the POD-based ROM. The $k_{\text{eff}}$ value for the AE-based ROM, seen in Figure~\ref{fig:2dae_keff_unseen}, has a final value close to that of the high-fidelity model, although convergence is not as smooth for the other models. 

\begin{figure}[!htb]
\centering
\begin{minipage}{.45\textwidth}
  \centering
  \includegraphics[scale=0.5]{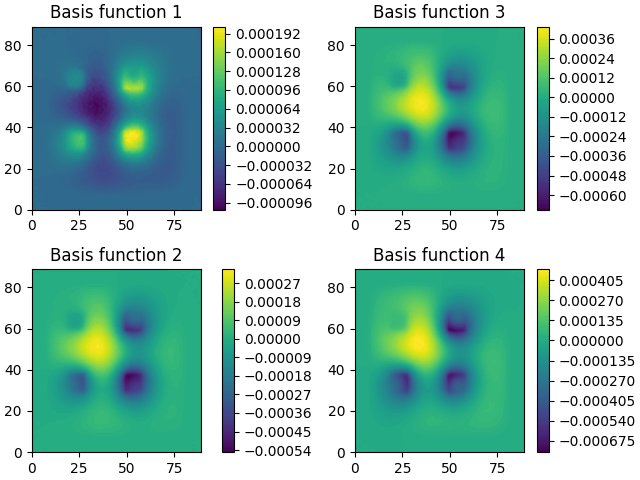}
  \subcaption{For the seen problem with parameters $r_1=0.9782$, $r_2=0.9891$, $r_3=0.7006$ and $r_4=0.8316$.}
  \label{fig:ae_4_basis_functions_seen}
\end{minipage}%
\hfill
\begin{minipage}{.45\textwidth}
  \centering
  \includegraphics[scale=0.5]{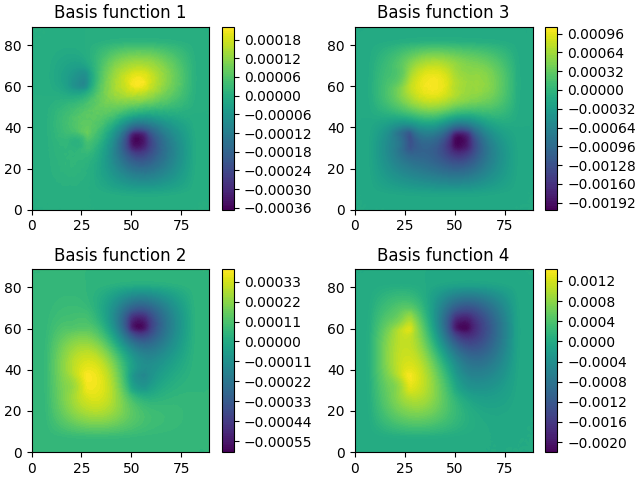}
  \subcaption{For the unseen problem with parameters $r_1=0.9452$, $r_2=0.8647$, $r_3=0.9996$ and $r_4=0.9776$.}
  \label{fig:ae_4_basis_functions_unseen}
\end{minipage}
\caption{The four basis functions from the linearised $\bm{C}$ matrix generated by the autoencoder at the point of convergence.}
\label{2daebasis}
\end{figure}

Figure~\ref{2daebasis} shows four basis functions of the nonlinear map, linearised at convergence, for both the seen and unseen sets of parameter values. These basis functions are columns of the $\bm{C}$ matrix. As for the POD-based model, it can be seen the basis functions capture variation associated with the control-rod regions. 

\begin{figure}[!htb]
\centering
\begin{minipage}{.45\textwidth}
  \centering
  \includegraphics[scale=0.5]{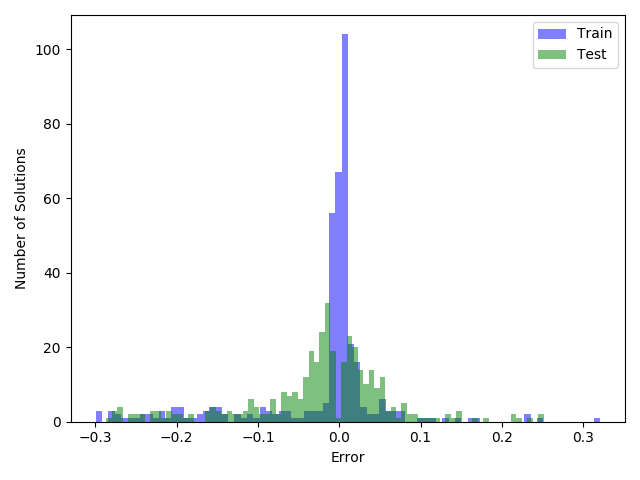}
  \subcaption{Histogram of the error in the flux profile, after removing outliers (5\% of the values).}
  \label{fig:2dAE_flux}
\end{minipage}%
\hfill
\begin{minipage}{.45\textwidth}
  \centering
  \includegraphics[scale=0.5]{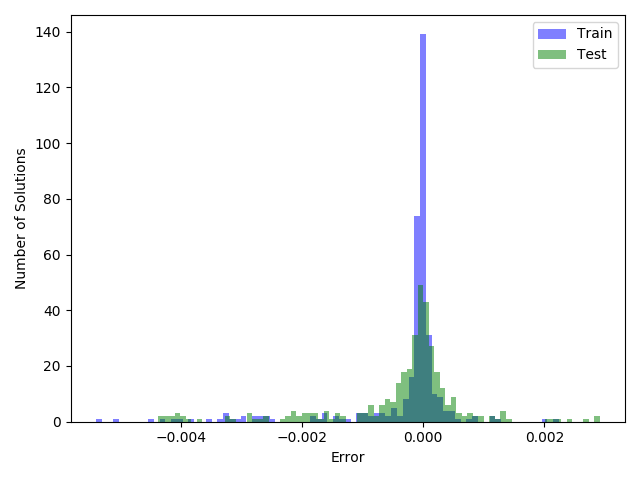}
  \subcaption{Histogram of the error in $k_{\text{eff}}$, after removing some outliers (20\% of values).}
  \label{fig:2dae_keff}
\end{minipage}
\caption{Error in the seen (blue) and the unseen (green) results for the autoencoder-based ROM applied to the 2D test case.}
\label{2daeerror}
\end{figure}

Figure~\ref{2daeerror} shows, despite the over-fitting of the network, that the training and test data both have similar ranges of errors as calculated by Equation~\eqref{error} for the error in scalar flux and Equation~\eqref{keff_error} for the error in $k_{\text{eff}}$. In comparison with results from the POD-based ROM, this model has a lower standard deviation associated with the errors, but a much longer tail. For this reason 5\% of the results have been removed from Figure~\ref{fig:2dAE_flux} to focus the graph on the more frequently-occurring errors. The original range of the errors was between -0.6 and +0.32. Figure~\ref{fig:2dae_keff} shows errors $k_{\text{eff}}$ from 80\% of the results,  where the original range of errors was between -0.04 and +0.003. The majority of the errors are lower than for the POD-based ROM, and the standard deviation is small compared to the results of the POD-based ROM.  

\subsubsection{Singular Value Decomposition and Autoencoder approach}\label{2d_svd_ae} 
As an alternative to having a large difference between  the number of neurons in the input and the first hidden layer, a combination of an SVD and an autoencoder is explored. First the SVD is applied to the snapshots, after which, 100 of the 400~possible basis functions are retained. The snapshots are projected onto the basis functions (see Equation~\eqref{alpha=rt}) and the resulting snapshot coefficients are used to train the autoencoder. The architecture is similar to that of the 1D problem, however, in order to be consistent with the other 2D reduced order models, the 100~variables are compressed to 4~latent variables. The autoencoder comprises of 14~fully connected layers with the number of neurons in each layer being: 
\begin{equation}
100 \rightarrow{} 100 \rightarrow{} 70 \rightarrow{} 50 \rightarrow{} 30 \rightarrow{} 16 \rightarrow{} 8 \rightarrow{} 4 \rightarrow{} 8 \rightarrow{} 16 \rightarrow{} 30 
\rightarrow{} 50 \rightarrow{} 70 
 \rightarrow{} 100 \rightarrow{} 100
    \nonumber
\end{equation}
The activation function for every layer was the exponential linear unit; the optimiser used was `Nadam' and the loss function was defined to be the mean squared reconstruction error. The snapshots were scaled between 0 and 1, and the network was trained using mini-batch gradient descent with a batch size of 50. The number of epochs was 30,000.

\begin{figure}[!htb]
\centering
\begin{minipage}{.45\textwidth}
  \centering
  \includegraphics[scale=0.5]{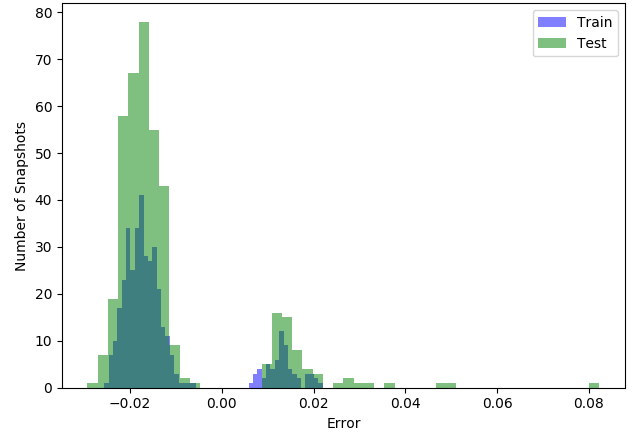}
  \subcaption{Reconstruction error in the flux profile for seen (blue) and unseen data (green).}
  \label{fig:svdae4pass_reconstruction}
\end{minipage}%
\hfill
\begin{minipage}{.45\textwidth}
  \centering
  \includegraphics[scale=0.5]{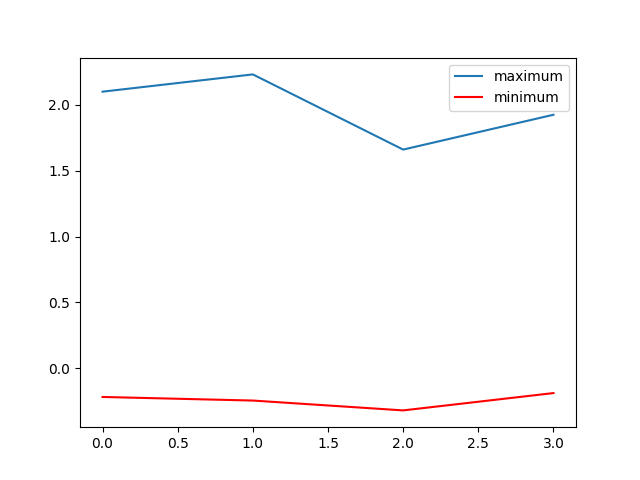}
  \subcaption{Minimum and maximum of the 4 latent variables generated from the SVD-autoencoder over the training data.}
  \label{fig:svdae4pass}
\end{minipage}
\caption{The representation error introduced to the representation of the high-fidelity model solutions when using the SVD-autoencoder, and the range of the latent variables.}
\label{svdae4passlat} 
\end{figure}

Figure~\ref{fig:svdae4pass_reconstruction} shows the error in the reconstruction when using the SVD-autoencoder for dimensionality reduction, calculated using Equation~\eqref{error_reconstruction_SVDAE}. 
The ranges of error for the training data and the test data are similar, which indicates that this network has been trained to a much better degree than the previous one, with no over-fitting occurring. Figure~\ref{fig:svdae4pass} shows the minimum and maximum values attained by the latent variables over the training data, which shows the relative sizes and ranges of the latent variables.

\begin{figure}[!htb]
\centering
\begin{minipage}{.45\textwidth}
  \centering
  \includegraphics[scale=0.5]{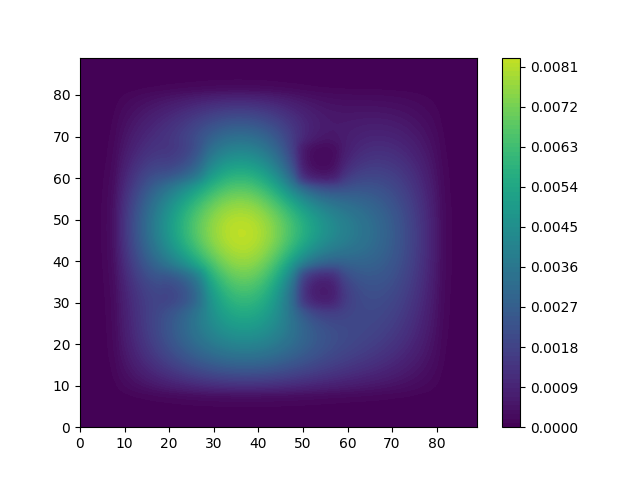}
  \subcaption{High-fidelity model.}
  \label{fig:dummyleft8}
\end{minipage}%
\hfill
\begin{minipage}{.45\textwidth}
  \centering
  \includegraphics[scale=0.5]{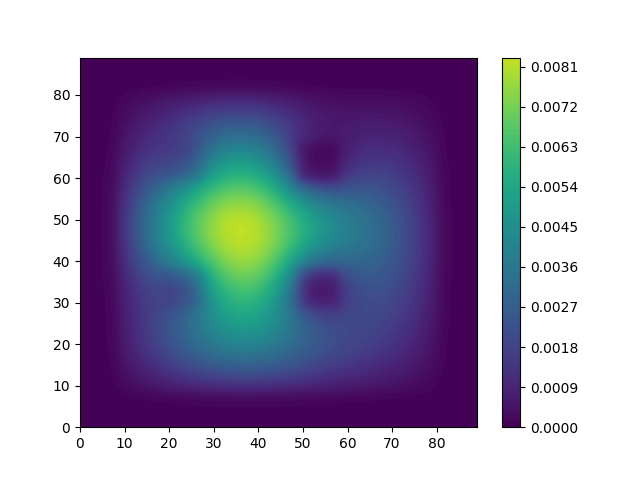}
  \subcaption{SVD-AE-based ROM.}
  \label{fig:dummyright8}
\end{minipage}
\caption{Scalar flux (\si{neutrons.cm^{-2}.s^{-1}}) across the reactor for a seen problem using the SVD-AE-based ROM with material properties $r_1=0.9782$, $r_2=0.9891$, $r_3=0.7006$ and $r_4=0.8316$ for the 2D test case. }
\label{2dsvdaeseen}
\end{figure}

\begin{figure}[!htb]
\centering
\begin{minipage}{.45\textwidth}
  \centering
  \includegraphics[scale=0.5]{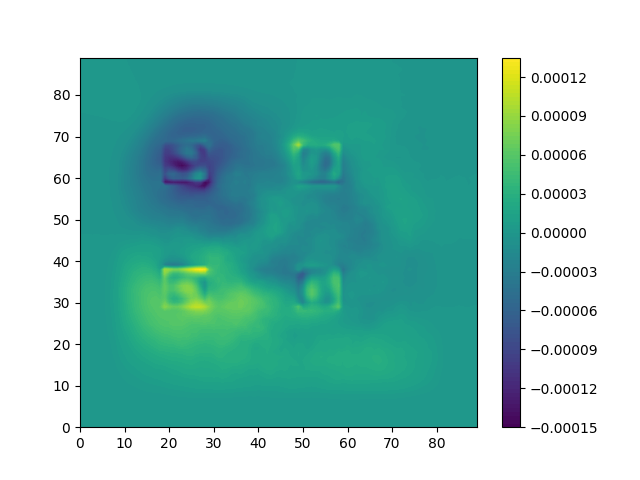}
  \subcaption{Pointwise error in scalar flux.}
  \label{fig:2dsvdae_diff_seen}
\end{minipage}%
\hfill
\begin{minipage}{.45\textwidth}
  \centering
  \includegraphics[scale=0.5]{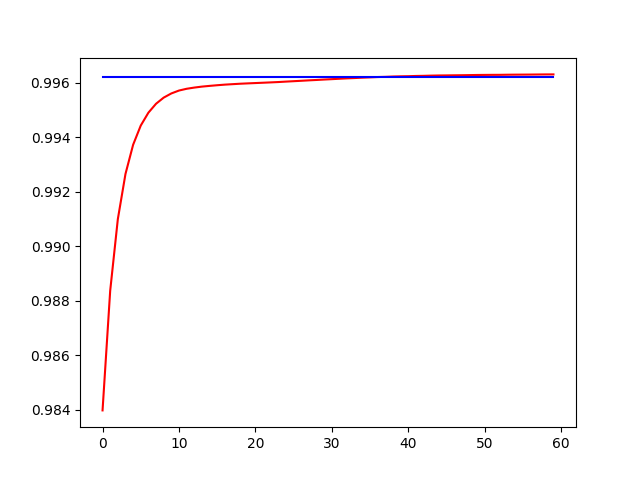}
  \subcaption{$k_{\text{eff}}$ vs iteration number, convergence of POD-based ROM (blue) to the converged HFM value (red).}
  \label{fig:2dsvdae_keff_seen}
\end{minipage}
\caption{Pointwise error in scalar flux and convergence of $k_{\text{eff}}$ for the SVD-AE-based ROM with material properties $r_1=0.9782$, $r_2=0.9891$, $r_3=0.7006$ and $r_4=0.8316$.}
\label{svdae_keffconv_seen}
\end{figure}

Figure~\ref{2dsvdaeseen} shows the scalar flux profiles of the high-fidelity model and the SVD-AE-based ROM for a seen set of parameter values. The agreement is good, and much better than the previous POD-based and AE-based ROMs. The pointwise error between the high-fidelity model and the AE-based ROM seen in Figure~\ref{fig:2dsvdae_diff_seen} is lower than for the two previous models and the error is concentrated around the control-rod regions. The $k_{\text{eff}}$ value of the SVD-AE-based model, seen in Figure~\ref{fig:2dsvdae_keff_seen}, converges to a value extrememly close to that of the high-fidelity model. 

\begin{figure}[!htb]
\centering
\begin{minipage}{.45\textwidth}
  \centering
  \includegraphics[scale=0.5]{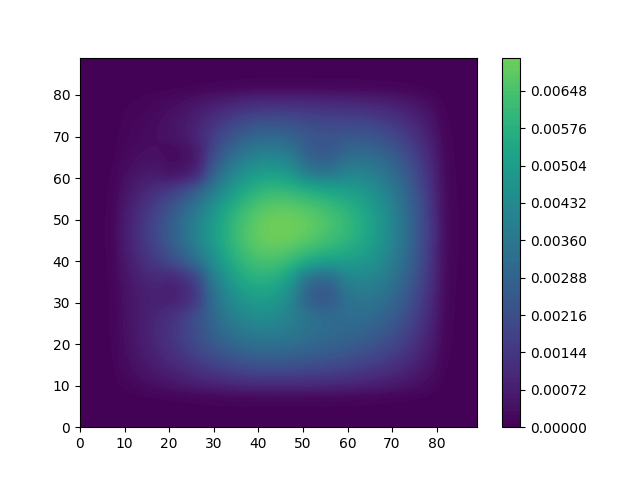}
  \subcaption{High-fidelity model.}
  \label{fig:dummyleft9}
\end{minipage}%
\hfill
\begin{minipage}{.45\textwidth}
  \centering
  \includegraphics[scale=0.5]{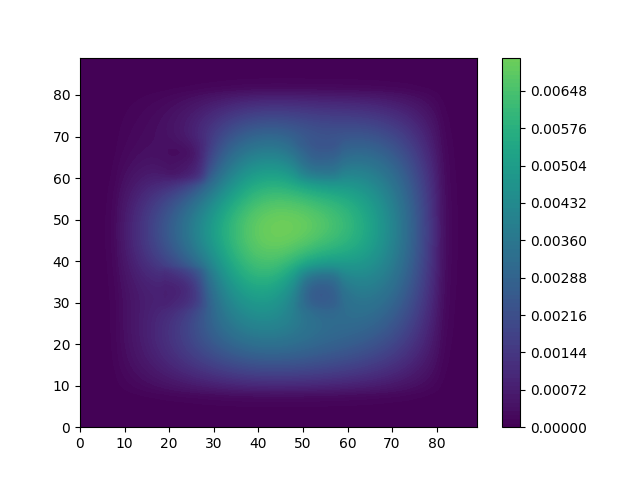}
  \subcaption{SVD-AE-based ROM}
  \label{fig:dummyright9}
\end{minipage}
\caption{Scalar flux (\si{neutrons.cm^{-2}.s^{-1}}) across the reactor for an unseen problem using the SVD-AE-based ROM with material properties $r_1=0.9452$, $r_2=0.8647$, $r_3=0.9996$ and $r_4=0.9776$.}
\label{2dsvdaeunseen}
\end{figure}

\begin{figure}[!htb]
\centering
\begin{minipage}{.45\textwidth}
  \centering
  \includegraphics[scale=0.5]{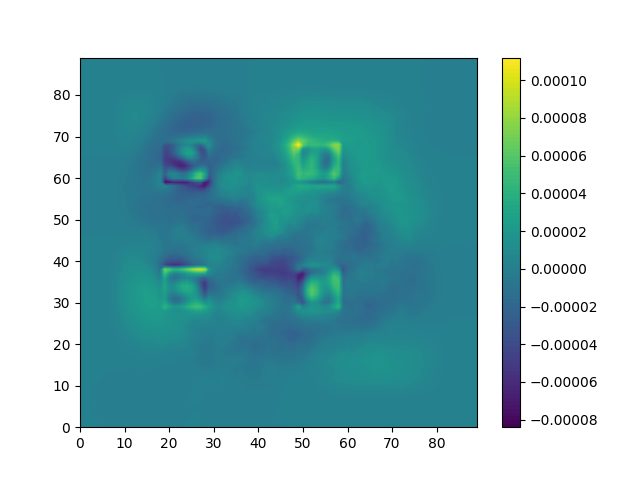}
  \subcaption{Pointwise error in scalar flux across the reactor.}
  \label{fig:2dsvdae_diff_unseen}
\end{minipage}%
\hfill
\begin{minipage}{.45\textwidth}
  \centering
  \includegraphics[scale=0.5]{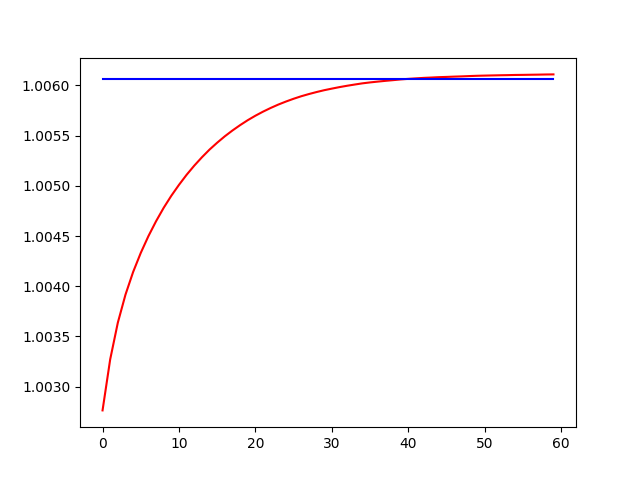}
  \subcaption{$k_{\text{eff}}$ vs iteration number, convergence of SVD-AE-based ROM (blue) to the converged HFM value (red).}
  \label{fig:2dsvdae_keff_unseen}
\end{minipage}
\caption{Pointwise error in scalar flux and convergence of $k_{\text{eff}}$ for an unseen problem using the SVD-AE-based ROM with material properties $r_1=0.9452$, $r_2=0.8647$, $r_3=0.9996$ and $r_4=0.9776$.}
\label{svdae_conv_keff_unseen}
\end{figure}

Figure~\ref{2dsvdaeunseen} shows results for an unseen set of parameter values. The SVD-AE-based ROM yields a scalar flux profile very similar to that of the high-fidelity model, and the pointwise error, shown in Figure~\ref{fig:2dsvdae_diff_unseen}, is again lower than the POD-based and AE-based ROMs with the highest values concentrated around the control-rod regions.  The $k_{\text{eff}}$ value for this model converges to the value obtained by the high-fidelity model, seen in Figure~\ref{fig:2dsvdae_keff_unseen}. 

\begin{figure}[!htb]
\centering
\begin{minipage}{.45\textwidth}
  \centering
  \includegraphics[scale=0.5]{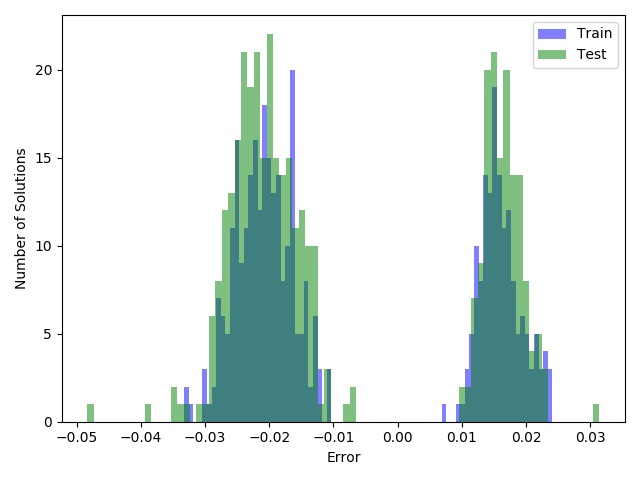}
  \subcaption{Histogram of the error in the flux profile.}
  \label{fig:dummyleft10}
\end{minipage}%
\hfill
\begin{minipage}{.45\textwidth}
  \centering
  \includegraphics[scale=0.5]{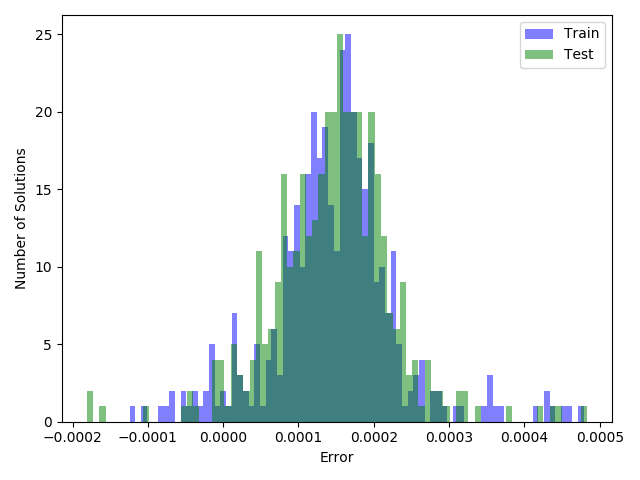}
  \subcaption{Histogram of the error in $k_{\text{eff}}$, after removing outliers (5\% of the data).}
  \label{fig:SVDaekeff}
\end{minipage}
\caption{Error in the seen (blue) and the unseen (green) results for the SVD-AE-based ROM applied to the 2D test case.}
\label{svdae4errn}
\end{figure}

The errors of the training and test data for the scalar flux and $k_{\text{eff}}$ can be seen in Figure~\ref{svdae4errn}. The errors are much lower  than those seen in both the POD-based ROM and the  AE-based ROM. Figure~\ref{fig:SVDaekeff} has had outliers removed (amounting to 5\% of the whole data set) to improve clarity of the plot. 

\begin{figure}[!htb]
\centering
\begin{minipage}{.45\textwidth}
  \centering
  \includegraphics[scale=0.5]{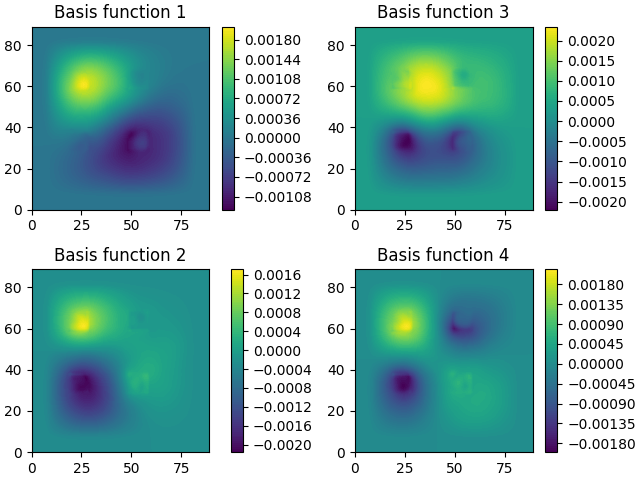}
  \subcaption{For the seen problem with parameters $r_1=0.9782$, $r_2=0.9891$, $r_3=0.7006$ and $r_4=0.8316$.}
  \label{fig:svdae_four_basis_functions_seen}
\end{minipage}%
\hfill
\begin{minipage}{.45\textwidth}
  \centering
  \includegraphics[scale=0.5]{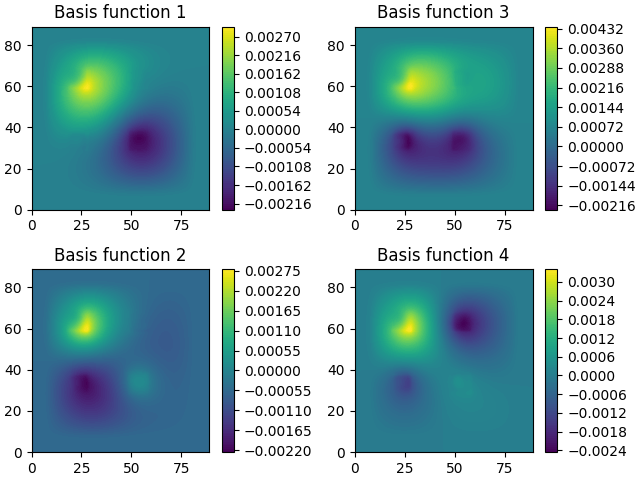}
  \subcaption{For the unseen problem with parameters $r_1=0.9452$, $r_2=0.8647$, $r_3=0.9996$ and $r_4=0.9776$.}
  \label{fig:svdae_four_basis_functions_unseen}
\end{minipage}
\caption{The four basis functions from the linearised $\bm{C}$ matrix generated by the SVD-autoencoder at the point of convergence.}
\label{2dsvdaebasis}
\end{figure}

Figure~\ref{2dsvdaebasis} shows four basis functions of the nonlinear map, linearised at the point of convergence, for both the seen and unseen sets of parameter values. These basis functions are columns of the $\bm{C}$ matrix. As for the POD-based ROM and the AE-based ROM, it can be seen the basis functions capture variation associated with the control-rod regions. 
 
Table~\ref{2derrors} shows the average maximum errors, defined in Equations~\eqref{average_maximum_error}, \eqref{average_maximum_error_keff} and~\eqref{average_maximum_error_recon},  for all three ROMs applied to the 2D problem. It can be seen that both autoencoder-based ROMs have lower compression errors than POD. However, the SVD-AE ROM performs better than both of AE-based ROM and POD-based ROM. Both AE-based ROMs generate better $k_{\text{eff}}$ values with the SVD-AE ROM being better than both the other models by an order of magnitude. The errors for the SVD-AE-based ROM are consistently lower than those of the other models.

\begin{table}[!htb]
\centering
\resizebox{\textwidth}{!}{%
 \begin{tabular}{l c c c c c c } 
 \hline\hline
         &  \multicolumn{2}{c}{POD} & \multicolumn{2}{c}{AE}& \multicolumn{2}{c}{SVD-AE}\\
 \hline\hline
         & Seen & Unseen & Seen & Unseen & Seen & Unseen  \\ 
 %\hline
Compression Error & \num{5.7397e-2} & \num{5.7756e-2} & \num{6.0167e-3} & \num{3.2919e-2} & \num{1.6785e-2} & \num{1.7734e-2}\\ 
 %\hline
 Flux Error & \num{6.4829e-2} & \num{6.5539e-2} & \num{5.9864e-2} & \num{7.8340e-2}  & \num{1.8995e-2} & \num{1.9395e-2}\\ 
 %\hline
 $k_{\text{eff}}$ Error & \num{2.1332e-3} & \num{1.9947e-3} & \num{1.6172e-3} & \num{1.8797e-3} & \num{1.5874e-4} & \num{1.7117e-4} \\
 \hline\hline
\end{tabular}
}
\caption{Average maximum errors for seen data (or training data) and unseen data (test data) for the POD-based ROM, AE-based ROM and SVD-AE-based ROM for the 2D test case. The number of compressed variables for all methods is four. }
\label{2derrors}
\end{table}

\clearpage
\section{Conclusions}\label{conc} 
By adopting an autoencoder for dimensionality reduction or compression, a novel projection-based reduced-order model for solving eigenvalue problems is developed. Two simple reactor problems are used as test cases to compare a number of autoencoder-based reduced order models with a reduced-order model (ROM) based on Proper Orthogonal Decomposition (POD). One of the autoencoder-based ROMs uses a novel combination of singular value decomposition (SVD) with an autoencoder.

POD and singular value decomposition (SVD) are often used to obtain basis functions of a low-dimensional or `reduced' space onto which the high-fidelity model is projected. This embedding of a high-dimensional space into a low-dimensional space is linear for POD, whereas, the embedding found by an autoencoder (AE) is, in general, nonlinear. This nonlinearity means that an autoencoder can provide a more accurate representation of data which could be important for larger, more complex problems.

For the 1D test case, when using a reduced space of dimension two, the AE-based reduced order model performs better than the POD-based model, with slightly lower errors in compression and in predicting the scalar flux profiles. When predicting $k_{\text{eff}}$, the AE-based ROM has an error that is one order of magnitude lower than the POD-based ROM. This supports the claim that each latent variable of the autoencoder can capture more information than each POD basis function.

%The rapid decay of the singular values in both test cases suggests that, with enough basis functions, the POD-based method will perform well. Indeed, for 10 basis functions, it outperforms the AE-based models in the 1D case. However, for the 1D case, when projecting onto a reduced space of dimension 2, the AE-based ROMs outperform the POD-based method. %
%This choice of dimension~2 was motivated by results from a ROM based on a variational autoencoder: such autoencoders are imbued with the ability to compress the data into the minimum number of latent variables. 
%For problems where many snapshots are used to capture the behaviour of the high-fidelity model, having an efficient representation of the low-dimensional space would be desirable. 

For the 2D test case, %in which the reduced space has a dimension of four, 
a standard autoencoder-based ROM demonstrated similar levels of accuracy to the POD-based ROM. However, the novel hybrid SVD-autoencoder is more accurate than both the other models, with errors that are an order of magnitude lower for $k_{\text{eff}}$. %The novel combination of SVD and an autoencoder \textcolor{red}{\textbf{helps}} with the accuracy of the method, as the length of the inputs to the autoencoder are much reduced.
Combining an SVD with an autoencoder reduces the length of the inputs to the autoencoder, which, in turn, reduces the complexity of the network, making overfitting less likely.

Although two relatively simple test cases are presented in this paper, our future work will involve increasingly complicated (i.e.~more realistic) assembly or reactor configurations which will benefit from the increased modelling capabilities that autoencoders provide. It can be seen that, for this eigenvalue problem, there is a clear benefit to using autoencoders over the standard POD approach. For minimal numbers of basis functions, autoencoders can be used to achieve a level of compression that cannot be matched in accuracy by POD. %The benefit to using autoencoders is highlighted by the fact that this is a linear problem where POD is known to perform well.

%\textcolor{red}{somewhat tangential: }The autoencoder's ability to recreate flux profiles is dependent on the training of the network and the errors for predicting the profiles are dependent on the errors associated with training. This means it is not unreasonable to assume that a better trained network has the potential to be more accurate than the networks seen here and, because of this, future work will focus on optimising the architecture of networks to improve results as well as applying the method to larger problems.

%The reason for this is that the autoencoder has a nonlinear representation of the solution variables and is thus more flexibly able to represent the solution than a linear representation from, for example, POD. 

\section*{Acknowledgements}
The authors would like to acknowledge the EPSRC Centre for Doctoral Training in Nuclear Energy (Imperial College, Cambridge University and the Open University,  EP/L015900/1) and the following EPSRC grants: MUFFINS (EP/P033180/1);  MAGIC (EP/N010221/1), INHALE (EP/T003189/1) and PREMIERE (EP/T000414/1).

\bibliographystyle{elsarticle-harv}%elsarticle-harv
\bibliography{ROM.bib} %,references_icferst

\end{document}